\documentclass[review,onefignum,onetabnum]{siamart251216}

\usepackage{lipsum}
\usepackage{amsfonts}
\usepackage{graphicx}
\usepackage{epstopdf}
\ifpdf
  \DeclareGraphicsExtensions{.eps,.pdf,.png,.jpg}
\else
  \DeclareGraphicsExtensions{.eps}
\fi

\newsiamremark{remark}{Remark}
\newsiamremark{example}{Example}
\newsiamremark{hypothesis}{Hypothesis}
\crefname{hypothesis}{Hypothesis}{Hypotheses}
\newsiamthm{claim}{Claim}

\headers{Orthogonal Constraint Pursuit}{ T. Li, D. Needell, and Q.W. Wang.}
\title{OCP: a greedy projection method for row-sparse representations of minimum-norm solutions\thanks{Corresponding author: Q.W. Wang (wqw@t.shu.edu.cn).
\funding{This work was funded by the National Natural Science Foundation of China [grant numbers 12401493 and 12371023].}}}
\author{Tao Li\thanks{School of Mathematics and Statistics, Hainan University, Haikou 570228, P.R. China
(\email{tli@hainanu.edu.cn}).}
 \and Deanna Needell\thanks{Department of Mathematics, University of California, Los Angeles 90095, USA  (\email{deanna@math.ucla.edu}).}
 \and Qing-Wen Wang\thanks{Department of Mathematics and Newtouch Center for Mathematics, and Collaborative Innovation Center for Marine Artificial Intelligence, Shanghai University, Shanghai 200444, P.R. China (\email{wqw@t.shu.edu.cn}).}}

\usepackage{amsopn}
\DeclareMathOperator{\diag}{diag}

\usepackage{algorithm}
\usepackage{algpseudocode}
\usepackage{booktabs,multirow,array}
\usepackage{amssymb}
\makeatletter
\newenvironment{breakablealgorithm}
  {\begin{center}
     \refstepcounter{algorithm}%
     \hrule height 0.8pt depth 0pt \kern 2pt
     \renewcommand{\caption}[2][\relax]{%
       {\raggedright\textbf{\ALG@name~\thealgorithm} ##2\par}%
       \ifx\relax##1\relax
         \addcontentsline{loa}{algorithm}{\protect\numberline{\thealgorithm}##2}%
       \else
         \addcontentsline{loa}{algorithm}{\protect\numberline{\thealgorithm}##1}%
       \fi
       \kern 2pt\hrule\kern 2pt}}
  {\kern 2pt\hrule\relax
   \end{center}}
\makeatother

\newcommand{\R}{\mathbb{R}}
\newcommand{\Sph}{\mathbb{S}}
\newcommand{\PP}{\mathbb{P}}
\newcommand{\EE}{\mathbb{E}}
\newcommand{\Range}{\mathcal{R}}
\newcommand{\Null}{\mathcal{N}}
\newcommand{\nnz}{\operatorname{nnz}}
\newcommand{\rank}{\operatorname{rank}}
\newcommand{\supp}{\operatorname{supp}}
\newcommand{\e}{\mathbf{e}}
\DeclareMathOperator*{\argmin}{arg\,min}
\DeclareMathOperator*{\argmax}{arg\,max}

\newif\ifrevisionmarks
\revisionmarkstrue
\ifdefined\mathcolor
\DeclareRobustCommand{\rev}[1]{\ifrevisionmarks\ifmmode\mathcolor{red}{#1}\else{\color{red}#1}\fi\else#1\fi}
\else 
\DeclareRobustCommand{\rev}[1]{\ifrevisionmarks{\color{red}#1}\else#1\fi}
\fi

\definecolor{revgreen}{rgb}{0,0.5,0}
\ifdefined\mathcolor
\DeclareRobustCommand{\revg}[1]{\ifrevisionmarks\ifmmode\mathcolor{revgreen}{#1}\else{\color{revgreen}#1}\fi\else#1\fi}
\else
\DeclareRobustCommand{\revg}[1]{\ifrevisionmarks{\color{revgreen}#1}\else#1\fi}
\fi

\begin{document}

\maketitle

\begin{abstract}
In this paper, we propose a greedy projection method, termed orthogonal constraint pursuit (OCP), for identifying a small set of rows of a consistent linear system $Ax=b$ whose span contains the minimum-norm solution $x^*$, along with the coefficients of that representation. The method successively activates the most violated constraint, and then projects the iterate onto the intersection of all activated constraints. We prove that OCP is equivalent to OMP applied to the dictionary $A^T$ and the signal $x^*$, and that it terminates in at most $\rank(A)$ iterations. When $x^*$ is a linear combination of $s$ rows of a matrix whose rows are independent and uniformly distributed on the unit sphere, OCP identifies those rows in exactly $s$ iterations with probability $1-O(m^{-1})$, provided that $n\ge Cs\ln m$ for a sufficiently large constant $C$. We also derive a variant of OCP with constraint removal for hard-margin support vector machines with regularized bias, which returns the support vectors, along with the multipliers that satisfy the KKT conditions. Numerical results on both real and synthetic problems demonstrate the identification of representing rows and the regime in which it is competitive with Krylov subspace methods.

\end{abstract}



\section{Introduction}\label{sec:intro}
In this paper, we consider the minimum-norm solution of a highly overdetermined and consistent linear system
\begin{equation}\label{eq-1}
Ax=b,
\end{equation}
where $A\in\R^{m\times n}$ and $b\in\R^m$ are given, with $m\gg n$. More precisely, we focus on representing this solution $x^*$ as a linear combination of a few rows of $A$. In the data analysis community, the rows of $A$ are commonly referred to as samples \cite{SeberLee}. Each nonzero row, together with the corresponding entry of $b$, defines a linear equation, whose solution set is a hyperplane in  $\R^n$. The well-known Kaczmarz-type methods, discussed below, project the iterate onto these hyperplanes. In our setting, we call these equations constraints, and identify each constraint with its hyperplane. Indeed, the constraints reduce the degrees of freedom in the solution space, each of them imposing a linear condition on the unknowns. A data set is a collection of such constraints, and their coefficient vectors form the rows of $A$. Naturally, finding such a representation is equivalent to extracting a few constraints that determine $x^*$ by themselves. We make this precise in the following definition.

\begin{definition}\label{def-1}
For $\Lambda\subset\{1,\dots,m\}$, let $A_\Lambda$ be the submatrix of $A$ formed by the rows indexed by $\Lambda$. We call  $\Lambda$ a \emph{representing set} for $x^*$ if $x^*\in\Range(A_\Lambda^T)$, i.e., $x^*=A_\Lambda^Ty_\Lambda$ for some $y_\Lambda\in\R^{|\Lambda|}$. A representing set is \emph{minimal} if it has the smallest cardinality among all representing sets for $x^*$, that is, if $|\Lambda|=\min\{\|y\|_0\ |\ x^*=A^Ty, y\in\R^m\}$.
\end{definition}

From Definition \ref{def-1}, the consistent subsystem $A_\Lambda x = b_\Lambda$ has a unique minimum-norm solution in $\Range(A_\Lambda^T)$, with the corresponding subvector of $b$ being $b_\Lambda$. 
Therefore, $\Lambda$ is a representing set if and only if this minimum-norm solution coincides with $x^*$; in this sense, the constraints in $\Lambda$ determine $x^*$ by themselves. In \cite{BDM13}, the selected
rows are weighted so that they approximate the residual norm $\|b-Ax\|_2$ uniformly in $x$, whereas a representing set reproduces
the single vector $x^*$ exactly. Finding a minimal representing set is equivalent to the following constrained optimization problem
\begin{equation}\label{eq-2}
\min_{y\in\R^m}\|y\|_0,\quad \text{s.t.}\quad AA^Ty=b,
\end{equation}
where $\|\cdot\|_0$ denotes the $\ell_0$-norm, which counts the nonzero entries of $y$. Remarkably, $AA^Ty=b$ holds if and only if  $A^Ty=x^*$, since $A^Ty\in\Range(A^T)$ and the minimum-norm solution is unique. Thus, \eqref {eq-2} can be viewed as finding a sparse representation of $x^*$ in the dictionary $A^T$. Conversely, the problem $\min\|y\|_0$ subject to $\Phi y=v$ with $v\in\Range(\Phi)$ is of the form \eqref{eq-2}, with $A=\Phi^T$ and $b=\Phi^Tv$. This shows that \eqref{eq-2} is NP-hard in general \cite{Natarajan95}. Moreover, if the minimum in \eqref{eq-2} is $s$,  then a minimal representing set $\Lambda_s$ consists of $s$ constraints, and $s\le\rank(A)$. Instead of relaxing the $\ell_0$-norm to the $\ell_1$-norm and solving the resulting convex problem, e.g., basis pursuit \cite{CandesTao05, RudelsonVershynin05}, we wish to formulate a greedy method, similar to orthogonal matching pursuit (OMP) \cite{Tropp04}, for finding a minimal representing set and its corresponding representation. This naturally gives rise to the following three questions:
\begin{enumerate}
\item Why is such a representation necessary?
\item How many constraints are required to reconstruct the solution?
\item What methods could implement the task efficiently?
\end{enumerate}

For the first question, both the pseudo-inverse and the linear combination of the $s$ rows of a minimal representing set yield the same minimum-norm solution $x^*$ of \eqref{eq-1}, but they generally correspond to different coefficient vectors in $\R^m$. The former, given by
$A^\dagger b=A^T(AA^T)^\dagger b$, corresponds to the dense coefficient vector $(AA^T)^\dagger b$. It is the minimum $\ell_2$-norm solution to $ A^Ty=x^*$, and distributes the weight over all rows. The latter corresponds to an $\ell_0$-minimizer, whose support is a minimal representing set. The two representations are algebraically identical if and only if an $\ell_0$-minimizer $y$ of \eqref{eq-2} lies in the column space of $A$, which is not generally the case. Once a minimal representing set $ \Lambda_s$ has been found, the solution $x^*$ is determined by the $s$ indices and their coefficients. It is recovered by one product with $A_{\Lambda_s}^T$, rather than a product with $A^T$.
The selection of a few rows that carry the information of a linear model arises in medical imaging \cite{LustigDonohoSantos}, data compression \cite{BDM13}, optimal sensor placement \cite{BruntonBruntonProctor}, and support vector machine (SVM) \cite{TsangKwokCristianini}. For instance, let $\mathcal D=\{(p_i,q_i)\}_{i=1}^m$ be a linearly separable training data set, in which $p_i\in\R^n$ are the feature vectors and $q_i\in\{-1,1\}$ are the binary labels. Assume that  $a_i=[q_ip_i^T,q_i]^T\in\R^{n+1}$ are the label-augmented feature vectors, and $A\in\R^{m\times(n+1)}$ is the coefficient matrix whose $i$th row is $a_i^T$. The hard-margin SVM \cite{Brooks11} seeks an optimal separating hyperplane, parameterized by a normal vector $w\in\R^n$ and a bias $\beta\in\R$, which maximizes the margin $1/\|w\|_2$ while classifying all samples correctly, i.e.,
\[
\min_{w\in\R^n,\beta\in\R}\ \frac{1}{2}\|w\|_2^2\quad\mbox{s.t.}\quad q_i(w^Tp_i+\beta)\ge1,\quad \forall i=1,\dots,m.
\]
We adopt the augmented formulation with $x=[w^T,\beta]^T\in\R^{n+1}$,
in which the bias term is also regularized. Then, the objective is strongly convex on $\R^{n+1}$, and the equality constraint $\sum_iq_iy_i=0$ disappears from the dual problem \cite{MangasarianMusicant99}. The augmented problem is reformulated as
\begin{equation}\label{eq-3}
\min_{x\in\R^{n+1}}\ \frac{1}{2}\|x\|_2^2\quad\text{s.t.}\quad Ax\ge\e,
\end{equation}
where $\e\in\R^m$ is the vector of all ones. For separable data, the feasible set of \eqref{eq-3} is a nonempty closed polyhedron, and the strictly convex objective has a unique minimizer, denoted by $x^*$. Since the constraints are linear, the Karush–Kuhn–Tucker (KKT) conditions are necessary and sufficient for optimality. This means that the vector $x^*\in\R^{n+1}$ yields \eqref{eq-3} if and only if there exists a multiplier vector $y^*\in\R^m$ such that
\begin{equation}\label{eq:kkt}
x^*=A^Ty^*,\quad y^*\ge0,\quad Ax^*\ge\e,\quad y_i^*(a_i^Tx^*-1)=0\ \ (i=1,\dots,m).
\end{equation}
The samples with $a_i^Tx^*=1$ are the support vectors. According to complementary slackness, the $y_i^*$ vanishes for all non-support vectors. If $\Lambda=\supp(y^*)$, then $x^*=A_\Lambda^Ty_\Lambda^*$ and $A_\Lambda x^*=\e_\Lambda$. This shows that $x^*$ is the minimum-norm solution of the subsystem $A_\Lambda x=\e_\Lambda$ formed by these support vectors.
Consequently, the hard-margin SVM reduces to a row-selection problem over consistent subsystems. It suffices to find a set $\Lambda$ of rows such that the minimum-norm solution $x_\Lambda=A_\Lambda^Ty_\Lambda$ of $A_\Lambda x=\e_\Lambda$ satisfies two conditions: it is feasible, $Ax_\Lambda\ge\e$, and its coefficients are nonnegative. For inequality constraints, the greedy selection can take the one-sided constraint violation
\[
(1-a_i^Tx_k)_+=\max\{0,1-a_i^Tx_k\},
\]
and constraints whose coefficients become negative have to be removed again. The resulting variant of our method is given in section~\ref{sec:ineq}. It terminates once the KKT conditions \eqref{eq:kkt} hold, and these conditions certify the computed support.

For the second question, it is trivial to see that no fewer than $s$ constraints are required.
The difficulty is to identify them efficiently, which is the third question.
Over the past two decades, row-action solvers, particularly Kaczmarz-type methods, have gained significant popularity for computing $x^*$, owing to their low memory and ease of implementation. The methods approximate the solution by projecting the current iterate onto constraints selected one after another, either randomly or deterministically. A large literature has grown related to Kaczmarz-type methods, yielding substantial progress in both theoretical bounds and computational efficiency. 
In 2009, Strohmer and Vershynin \cite{StrohmerVershynin09} 
introduced the randomized Kaczmarz (RK) method, which 
selects the row used at each iteration with probability proportional to its squared Euclidean norm, rather than in a given order.  To improve its performance, Needell and Tropp \cite{NeedellTropp14} proposed a randomized block method; for other block variants, see \cite{Censor81,Necoara19,MiaoWu22} and the references therein. In 2018, Bai and Wu \cite{BaiWu18} developed the greedy randomized Kaczmarz (GRK) method, and 
proved a smaller convergence factor than RK.
Moreover, Du et al. \cite{DuHayamiZheng21} derived a Kaczmarz-type inner-iteration preconditioned flexible GMRES method for solving \eqref{eq-1}. To reduce the considerable cost of GRK, Jiang et al. \cite{JiangWuJiang23} devised a semi-randomized Kaczmarz method with simple random sampling. Derezi\'nski et al. \cite{DNR25} proposed an accelerated randomized block Kaczmarz method, which exploits outlying singular values in the input to attain fast, Krylov-style convergence. Recently, Jeong et al. \cite{JeongNeedellWu26} proposed a federated Kaczmarz method and a modified version, both of which are promising.

Note that Krylov subspace methods, such as CGNE \cite{Craig55}, CGLS \cite{Bjorck96}, and LSQR \cite{PaigeSaunders82}, terminate at $x^*$ in at most $\rank(A)$ iterations in exact arithmetic, but they cannot provide a sparse linear representation of $x^*$.
Specifically, starting with $x_0=0$, the CGNE iterates satisfy $x_k=A^Ty_k$ with $y_k=p_{k-1}(AA^T)\,b$, where $p_{k-1}$ is a polynomial of degree at most $k-1$. The vectors $y_k$ converge to $(AA^T)^\dagger b$, which is dense in general.  As a result, the iterates do not 
present a small set of rows that represents $x^*$. On the other hand, Kaczmarz-type methods use one row, or a block of rows, at a time and do not retain previously used rows, so their iterates do not provide such a set either. In applications such as support vector extraction in hard-margin SVM, sparse sensor placement, and data compression, however, what is required is the pair $(\Lambda,y_\Lambda)$ with $x^*=A_\Lambda^Ty_\Lambda$ and $|\Lambda|$ small, not merely the minimizer $x^*$. Indeed, $x^*$ is obtained by one product with $A_\Lambda^T$ once $(\Lambda,y_\Lambda)$ is known.

It is worth noting that the 
greedy column-action methods have been studied extensively in sparse approximation and compressed sensing. The orthogonal matching pursuit (OMP) method \cite{PRK93,DMA97} selects
the column (atom) of a dictionary that is most correlated with the current residual, and then projects the signal orthogonally onto the span of the selected atoms. In 2004, Tropp \cite{Tropp04} showed that OMP recovers the sparsest representation under the exact recovery condition, which holds in particular when the coherence of the dictionary is small enough. Moreover, Tropp and Gilbert \cite{TroppGilbert07} showed that OMP recovers an $s$-sparse vector in $\R^d$ from $O(s\ln d)$ random measurements with high probability; see \cite{NeedellVershynin09} for regularized OMP and \cite{NeedellTropp09} for CoSaMP. As noted above, \eqref{eq-2} is the sparse representation problem for 
$x^*$ in the dictionary $A^T$. Observe that OMP cannot be applied to it directly, since $x^*$ is not available in advance. Nevertheless, OMP uses the signal only through its correlations with the atoms, and for $x^*$ the data give these correlations, $Ax^*=b$. 
This leads to the method studied in this paper, which we call orthogonal constraint pursuit (OCP). It works in the row space of the coefficient matrix. The index of the largest entry of the residual in absolute value is appended to an index set at each iteration. The iterate is then projected onto the intersection of all selected constraints. We will illustrate in section~\ref{sec:ocp} that
OCP is OMP applied to $A^T$.
This correspondence allows us to transfer the analysis of OMP to the selection of rows of a linear system. Moreover, if $x^*$ is a combination of $s$ rows, a
natural question is whether OCP identifies these active rows in exactly $s$ iterations. We will give a probabilistic guarantee that it does, and a deterministic sufficient condition for exact recovery. Furthermore, OCP can be viewed as the maximal-residual relaxation method of Agmon \cite{Agmon54} and Motzkin and Schoenberg \cite{MotzkinSchoenberg54}, in which the projection onto the newest hyperplane is replaced with the projection onto the intersection of all selected constraints.
Unlike row-action methods, which are memoryless, it retains information from previously selected rows.

Actually, OCP is intended for systems whose rank $r = \rank(A)$ is far smaller than the number of rows $m$.
Every vector in the $r$-dimensional space $\Range(A^T)$, in particular the solution $x^*$, is a linear combination of at most $r$ rows. Therefore, it is determined by at most $r$ of the $m$ constraints. This number is even smaller when $x^*$ is a combination of only  $s < r$ rows, as in the random model of section~\ref{sec:recovery}. The proposed method finds a representing set with at most $r$ rows, together with its coefficients, in at most $r$ iterations. Its cost depends on the number of selected rows, and $m$ enters only through the products with $A$. The hard-margin SVM is an example. There, $r \le n+1 \ll m$, and the representing rows are support vectors. Identifying these vectors 
is then the primary goal, rather than merely finding 
$x^*$. The variant of OCP in section~\ref{sec:ineq} is devised for this purpose.
On the other hand, if $r$ is large and  $x^*$ admits no sparser row representation, the cost 
of OCP approaches that of a full QR factorization. Additionally, if only $x^*$ is required, Krylov subspace methods are commonly preferable. 

The rest of this paper is organized as follows.
In section~\ref{sec:ocp}, we formulate the orthogonal constraint pursuit method and establish its correspondence with OMP. We also show that the rows of the selected constraints are linearly independent, and OCP terminates at $x^*$ in at most $\rank(A)$ iterations. We present an implementation, built on updating a QR factorization of $A_\Lambda^T$, along with its cost and storage. For linear inequalities, we give a 
variant with constraint removal, which turns out to be the dual active-set method of Goldfarb and Idnani \cite{GoldfarbIdnani83} with the maximal-violation selection rule. In section~\ref{sec:recovery}, we present the recovery results for OCP. 
If the rows are independent and uniformly distributed on the unit sphere and $x^*$ is a linear combination of $s$ of them, with support and coefficients chosen independently of $A$, then OCP identifies these rows in exactly $s$ iterations with probability $1-O(m^{-1})$, provided $n\ge Cs\ln m$ for a sufficiently large constant $C$. In section~\ref{sec:numerics},
we present some numerical examples to illustrate the identification of representing rows and support vectors and demonstrate when OCP is competitive with Krylov subspace methods for computing $x^*$. Finally, we conclude the whole paper by giving some remarks.

\section{Orthogonal constraint pursuit method}\label{sec:ocp}

In this section, we propose the orthogonal constraint pursuit (OCP) method for finding a row-sparse representation of the minimum-norm solution of  \eqref{eq-1}.
To avoid computing the matrices $A_\Lambda A_\Lambda^T$ explicitly, we provide an implementation of OCP based on QR factorization updates. Moreover, we extend the method to linear inequalities. Before doing so, we briefly recall some helpful results which will be used in the sequel.

Let $\R^{m\times n}$ and $\R^n$ be the sets of all real matrices of size $m\times n$ and of all real $n$-dimensional vectors, respectively. The inner product of two $n$-dimensional real vectors $x$ and $y$ is defined as $\langle x,y\rangle=\sum_{i=1}^nx_iy_i$, and the induced $\ell_2$-norm of $x$ is $\|x\|_2=\sqrt{\sum_{i=1}^nx_i^2}$. The $\ell_\infty$-norm of $x$ is $\|x\|_{\infty}=\max\{|x|_1,\cdots,|x|_n\}$.
Denote by $A^T$, $A^{-1}$, and $A^\dagger$ the transpose, inverse, and pseudo-inverse of $A$, respectively. The rows of a matrix $A\in\R^{m\times n}$ are denoted by $a_1^T,a_2^T,\dots,a_m^T$, its null space by $\Null(A)$, and the range space of $A^T$ by $\Range(A^T)$. The latter is spanned by the column vectors $a_i$, $i=1,\dots,m$. The $i$th constraint equation thus reads $a_i^Tx=b_i$, and we identify this constraint with the hyperplane $\{x\in\R^n\mid a_i^Tx=b_i\}$, whose normal vector is $a_i$. The feasible set, which consists of the solutions $\tilde x$ of \eqref{eq-1}, is of the form $\Omega=\{\tilde x\in\R^n\mid A\tilde x=b\}$. Recall from \cite{Saad03} that $\R^n$ admits the orthogonal direct sum decomposition
\[
\R^n=\Range(A^T)\oplus\Null(A),
\]
which leads to the following result.

\begin{proposition}\label{prop:decomp}
The intersection of $\Omega$ and $\Range(A^T)$, denoted by $\Omega\cap\Range(A^T)$,  consists of a single point, denoted by $x^*$, which is the unique minimum-norm solution of \eqref{eq-1}, and $\Omega=x^*+\Null(A)$. Every $\tilde x\in\Omega$ satisfies $\|\tilde x\|_2^2=\|x^*\|_2^2+\|\tilde x-x^*\|_2^2$.
\end{proposition}

\begin{proof}
Let $\tilde x\in\Omega$, which is nonempty by consistency, and let $\tilde x=u+z$ with $u\in\Range(A^T)$ and $z\in\Null(A)$.
 Then $Au=A\tilde x=b$, so that $u\in\Omega\cap\Range(A^T)$. If $u'$ also belongs to this set, then $u-u'\in\Range(A^T)\cap\Null(A)=\{0\}$. Thus,  $\Omega\cap\Range(A^T)=\{u\}$, and we let $x^*=u$. For every $\tilde x\in\Omega$, it follows that $\tilde x-x^*\in\Null(A)$, which is orthogonal to $x^*$. This gives the remaining assertions.
\end{proof}

\begin{definition}\label{def:sub}
Given an index set $\Lambda\subset\{1,2,\dots,m\}$, let $A_\Lambda\in\R^{|\Lambda|\times n}$ be the submatrix defined in Definition~\ref{def-1}, and $b_\Lambda\in\R^{|\Lambda|}$ be the corresponding subvector. The localized constraint subsystem is of the form  $A_\Lambda x=b_\Lambda$.
\end{definition}

Without loss of generality,
we assume that the rows of $A$ are $\ell_2$-normalized, i.e., $\|a_i\|_2=1$ for all $i=1,\dots,m$. This is no restriction. A zero row $a_i^T$ has $b_i=0$, since  \eqref{eq-1} is consistent, and can be deleted. Scaling an equation $a_i^Tx=b_i$ by a nonzero factor changes neither $\Omega$ nor $x^*$, nor the representing sets.
It only rescales the coefficients of the representations. In this setting, the algebraic residual of a constraint equals, up to sign, the distance to its hyperplane, so that the greedy selection becomes scale-invariant. From Proposition~\ref{prop:decomp}, the $x^*$, lying in  $\Range(A^T)$, is surely a linear combination of the rows of $A$. Therefore, for
 \eqref{eq-1}, it is natural to formulate the iterative scheme
\begin{equation}\label{eq:iter}
x_{k+1}=x_k+\lambda_ka_{i_{k+1}},
\end{equation}
where $\lambda_k$ is the step length, and $a_{i_{k+1}}$ is the normal vector of the constraint selected at this step. Let $e_k=x^*-x_k$ be the error
vector at the $k$th iterate, and $r_k=b-Ax_k=Ae_k$  be its corresponding residual vector. Then $r_k(i)=\langle a_i,e_k\rangle$ is the $i$th entry of $r_k$. Now we want to find two optimal
parameters $\lambda_k$ and $i_{k+1}$ yielding
\begin{align*}
\min\|e_{k+1}\|_2^2&=\min_{\lambda_k,i_{k+1}}\|e_k-\lambda_ka_{i_{k+1}}\|_2^2\\
&=\min_{\lambda_k,i_{k+1}}\bigl[\|e_k\|_2^2-2\lambda_k\langle e_k,a_{i_{k+1}}\rangle+\lambda_k^2\bigr].
\end{align*}
For a fixed index $i$, the above quadratic function of $\lambda_k$ attains its minimum $\|e_k\|_2^2-r_k(i)^2$ at $\lambda_k=\langle e_k,a_i\rangle=r_k(i)$. This gives that $\lambda_k=r_k(i_{k+1})$, and the optimal index is
\begin{equation}\label{eq:select}
i_{k+1}=\argmax_{1\le i\le m}|r_k(i)|=\argmax_{1\le i\le m}|b_i-a_i^Tx_k|=\argmax_{1\le i\le m}|\langle x^*-x_k,a_i\rangle|.
\end{equation}
With this relation, we can obtain that \eqref{eq:iter} is the maximal-residual relaxation method of Agmon \cite{Agmon54} and Motzkin and Schoenberg \cite{MotzkinSchoenberg54},
often called the greedy Kaczmarz method; for its analysis, see \cite{DeLoeraHaddockNeedell17,HaddockNeedell19}.

As stated earlier, the randomized Kaczmarz-type solvers approximate the minimum-norm solution, but they neither retain the previously selected rows nor produce a row-sparse representation. Our idea is to select constraints by \eqref{eq:select}, as in the
 maximal-residual method, but keeping all of them.
  At each iteration, we choose the constraint with the largest residual entry in magnitude, that is, the constraint whose row is most correlated with the current error $e_k$, and append its index to the active set, which contains the indices of the selected constraints so far. Then we project the iterate onto the intersection of all selected constraints. The iteration terminates when the residual vanishes, that is, when the active set is a representing set for $x^*$. When $x^*$ has a representing set of $s$ rows, the proposed method can find it in $s$ iterations. The conditions under which this happens are given in section~\ref{sec:recovery}. This gives the following method, called orthogonal constraint pursuit.
  \vspace{4mm}
\begin{breakablealgorithm}
\caption{Orthogonal Constraint Pursuit}\label{alg:ocp}
\begin{algorithmic}[1]
\Require Input $A$, $b$, the tolerance $\varepsilon\ge0$, and the maximum number of iterations $k_{\max}\ge\rank(A)$.
\State \textbf{Initialization:} Set the initial guess $x_0=0$, the active set $\Lambda_0=\emptyset$, the residual $r_0=b$, and the iteration counter $k=0$.
\While{$\|r_k\|_2>\varepsilon$ \textbf{and} $|\Lambda_k|<k_{\max}$}
\State \textbf{Step A (Greedy Constraint Selection):} Find the index $i_{k+1}$ of the constraint with maximum violation
\[
i_{k+1}=\argmax_{1\le i\le m}|r_k(i)|=\argmax_{1\le i\le m}|b_i-a_i^Tx_k|.
\]
If the maximum is attained at several indices, break the tie deterministically.
\State \textbf{Step B (Support Augmentation):} Update the active set and the matrix of the selected rows: $\Lambda_{k+1}=\Lambda_k\cup\{i_{k+1}\}$ and $A_{\Lambda_{k+1}}=\begin{bmatrix}A_{\Lambda_k}\\ a_{i_{k+1}}^T\end{bmatrix}$, where $A_{\Lambda_0}$ is the empty matrix.
\State \textbf{Step C (Subspace Orthogonal Projection):} Compute the minimum-norm solution of the subsystem formed by the active constraints
\begin{equation}\label{eq:stepC}
x_{k+1}=\argmin_{x\in\R^n}\|x\|_2\quad\text{s.t.}\quad A_{\Lambda_{k+1}}x=b_{\Lambda_{k+1}},
\end{equation}
where $b_{\Lambda_{k+1}}$ is the real vector of size $k+1$ whose $j$th component is the $i_j$th component of the right-hand side vector $b$, $j=1,\dots,k+1$, or equivalently,
\begin{equation}\label{eq:stepC2}
x_{k+1}=A_{\Lambda_{k+1}}^Ty_{k+1},\qquad y_{k+1}=\argmin_{y\in\R^{k+1}}\|A_{\Lambda_{k+1}}^Ty-x^*\|_2.
\end{equation}
The equivalence of \eqref{eq:stepC} and \eqref{eq:stepC2} is shown in Proposition~\ref{prop:omp}, and $y_{k+1}$ is unique by Lemma~\ref{lem:indep}.
\State \textbf{Step D (Residual Update):} Update the residual vector: $r_{k+1}=b-Ax_{k+1}$.
\State Increment the iteration counter: $k\leftarrow k+1$.
\EndWhile
\Ensure Output the index set $\Lambda_k$, the sparse vector $y\in\R^m$ supported on $\Lambda_k$ whose subvector on $\Lambda_k$ is $y_k=(A_{\Lambda_k}A_{\Lambda_k}^T)^{-1}b_{\Lambda_k}$, and $x_k=A_{\Lambda_k}^Ty_k$.
\end{algorithmic}
\end{breakablealgorithm}
 \vspace{4mm}

Observe from Algorithm~\ref{alg:ocp} that, the bound $|\Lambda_k|<k_{\max}$ is not necessary in exact arithmetic with   $\varepsilon=0$, since the residual vanishes after at most $\rank(A)$ iterations.
It guards against an endless loop in floating-point arithmetic, and $k_{\max}$ can be any upper bound for $\rank(A)$, such as $\min\{m,n\}$.

The following proposition indicates that Algorithm~\ref{alg:ocp} is OMP applied to the dictionary $A^T$ and the signal $x^*$. Recall that OMP starts from  $v_0=0$ and $\Lambda_0^{\rm OMP}=\emptyset$.  It selects the atom most correlated with the residual $x^*-v_k$ per iteration, and sets $v_{k+1}$ to be the orthogonal projection of 
 $x^*$ onto the span of the selected atoms.

\begin{proposition}\label{prop:omp}
Let $x^*$ be the unique minimum-norm solution of 
 \eqref{eq-1}, and let $\{\Lambda_k\}$, $\{x_k\}$ be generated by Algorithm~\ref{alg:ocp}. If $\{\Lambda_k^{\rm OMP}\}$ and $\{v_k\}$ are generated by OMP applied to the dictionary $\Phi=A^T=[a_1,\dots,a_m]$ and the signal $x^*$, with the same tie-breaking rule, then
\[
\Lambda_k=\Lambda_k^{\rm OMP},\qquad x_k=v_k=P_{\Range(A_{\Lambda_k}^T)}x^*,\qquad r_k=A(x^*-x_k).
\]
\end{proposition}

\begin{proof}
Since $r_k=b-Ax_k=A(x^*-x_k)$, it follows that $r_k(i)=\langle a_i,x^*-x_k\rangle$ is the correlation of the $i$th atom with the OMP error $x^*-x_k$. This indicates that Step~A and the selection step of OMP follow the same index whenever $x_k=v_k$. The feasible set of \eqref{eq:stepC} is $x^*+\Null(A_{\Lambda_{k+1}})$, and $\R^n=\Range(A_{\Lambda_{k+1}}^T)\oplus\Null(A_{\Lambda_{k+1}})$ is an orthogonal direct sum. Therefore, the point in this set with the minimal norm is $P_{\Range(A_{\Lambda_{k+1}}^T)}x^*$, i.e., the OMP iterate $v_{k+1}$. Starting with $x_0=v_0=0$, the proof can be completed by induction on $k$.
\end{proof}

As seen, Algorithm~\ref{alg:ocp} is OMP applied to $A^T$, and the only point to be checked is that the correlations required by OMP, $Ax^*$ for the projection and $A(x^*-x_k)$ for the selection, can be evaluated as $b$ and $b-Ax_k$ without knowledge of  $x^*$. The same observation underlies Batch-OMP \cite{RZE08} that carries out OMP using only $D^Tx^*$ and $D^TD$; for $D=A^T$, these are $b$ and $AA^T$.
We state the correspondence
explicitly because it allows the analysis of OMP to be applied to Algorithm~\ref{alg:ocp}.
From Proposition~\ref{prop:omp}, we have $x_k\in\Range(A_{\Lambda_k}^T)\subseteq\Range(A_{\Lambda_{k+1}}^T)$. Moreover, Step~C can be rewritten as $x_{k+1}=P_{L_{k+1}}(x_k)$, since the orthogonal projection onto the affine subspace $L_{k+1}=\{x\mid A_{\Lambda_{k+1}}x=b_{\Lambda_{k+1}}\}$ is of the form $P_{L_{k+1}}(z)=x_{k+1}+P_{\Null(A_{\Lambda_{k+1}})}z$. Replacing $L_{k+1}$ with the single constraint $\{x\mid a_{i_{k+1}}^Tx=b_{i_{k+1}}\}$ gives $x_{k+1}=x_k+r_k(i_{k+1})\,a_{i_{k+1}}$, the maximal-residual method. Consequently,  Algorithm~\ref{alg:ocp} can be viewed as a greedy row-space projection method, which incorporates the
maximal-residual rule with projection onto the intersection of all previously selected constraints.

In what follows, we present the basic properties of Algorithm~\ref{alg:ocp}. First, we show that the iterates satisfy all selected constraints, that their corresponding error decreases strictly, and that the iteration terminates at $x^*$. Proposition~\ref{prop:omp} leads to short proofs. The first one guarantees that the index $i_{k+1}$ chosen by the greedy selection is always new. 

\begin{lemma}\label{lem:active}
Let $x_k$ be the iterate generated by Algorithm~\ref{alg:ocp} at the $k$th iteration. Then
$x_k$ yields all previously activated constraints exactly, that is, $A_{\Lambda_k}x_k=b_{\Lambda_k}$. Hence, the localized residual vanishes on the support set, i.e., $r_k(j)=0$ for all $j\in\Lambda_k$.
\end{lemma}

\begin{proof}
From \eqref{eq:stepC}, its feasible set is a nonempty closed affine subspace.
The point with minimum norm lies in this set, so that $x_k$ satisfies the constraints of \eqref{eq:stepC} exactly, $A_{\Lambda_k}x_k=b_{\Lambda_k}$. This means that $r_k(j)=b_j-a_j^Tx_k=0$ for all $j\in\Lambda_k$.
\end{proof}

Putting this lemma with \eqref{eq:select} indicates that no index is selected twice. Indeed, before termination, we have $r_k\neq 0$, so that  $|r_k(i_{k+1})|>0$, whereas  $r_k(j)=0$ for $j\in\Lambda_k$. So, Algorithm~\ref{alg:ocp} terminates after at most
$m$ iterations. In fact, more is true.

\begin{lemma}\label{lem:indep}
Let $A_{\Lambda_k}$ be the active submatrix generated by Algorithm~\ref{alg:ocp} at the $k$th iteration. Then $A_{\Lambda_k}$ has full row rank for every $k$.
\end{lemma}

\begin{proof}
The proof is by induction on $k$. The result is trivial for $k=0$. Suppose it holds for $k$, and the loop is executed at step $k$. From Proposition ~\ref{prop:omp}, it follows that $(x^*-x_k)\perp\Range(A_{\Lambda_k}^T)$. If $a_i\in\Range(A_{\Lambda_k}^T)$, then
\[
r_k(i)=\langle a_i,x^*-x_k\rangle=0.
\]
Since the loop is executed, $r_k\neq0$, we have $|r_k(i_{k+1})|>0$.  Therefore, it follows that $a_{i_{k+1}}\notin\Range(A_{\Lambda_k}^T)$, i.e., the row  $a_{i_{k+1}}^T$ of the newly activated constraint is linearly independent of the rows of $A_{\Lambda_k}$, and $A_{\Lambda_{k+1}}$ has full row rank. 
\end{proof}

A consequence of the above lemma is the following proposition.

\begin{proposition}\label{prop:rowspace}
Let $A_{\Lambda_k}$ be the active submatrix defined in Lemma~\ref{lem:indep}, and let $\Range(A^T)$ be the range space of $A^T$. Then the iterative sequence $\{x_k\}$ generated by Algorithm~\ref{alg:ocp} remains in $\Range(A^T)$, i.e., $x_k\in\Range(A^T)$ for all $k\ge0$.
\end{proposition}

\begin{proof}
For $k=0$, the assertion is trivially true. For $k>0$, let $G_k=A_{\Lambda_k}A_{\Lambda_k}^T$, a real matrix of size $k\times k$. From Lemma~\ref{lem:indep}, we see that $G_k$ is invertible. This results in the minimum-norm solution, shown in Step~C, being expressed as $x_k=A_{\Lambda_k}^T(A_{\Lambda_k}A_{\Lambda_k}^T)^{-1}b_{\Lambda_k}$. With the vector $y_k=(A_{\Lambda_k}A_{\Lambda_k}^T)^{-1}b_{\Lambda_k}$ of dimension $|\Lambda_k|$, the $x_k$ is a linear combination of the columns of $A_{\Lambda_k}^T$, i.e.,
\begin{equation}\label{eq:xk}
x_k=A_{\Lambda_k}^Ty_k=\sum_{l=1}^{k}y_k(l)\,a_{i_l}\in\Range(A_{\Lambda_k}^T).
\end{equation}
Since $\Range(A_{\Lambda_k}^T)\subseteq\Range(A^T)$, it follows that $x_k\in\Range(A^T)$ for all iterations performed.
\end{proof}

Next, we demonstrate that  OCP terminates at the solution in at most $\rank(A)$ iterations, and its error norm decreases strictly.

\begin{lemma}\label{lem:monotone}
Let $x^*$ be the unique minimum-norm solution of  \eqref{eq-1}, and let $\{x_k\}$ be the iterative sequence generated by Algorithm~\ref{alg:ocp}. 
Then, for every $k$ at which the loop of Algorithm \ref{alg:ocp} is executed, its error norm decreases strictly
\begin{equation}\label{eq:monotone}
\|x^*-x_{k+1}\|_2<\|x^*-x_k\|_2.
\end{equation}
\end{lemma}

\begin{proof}
At the $k$th iteration, the error in Algorithm ~\ref{alg:ocp} is of the form
\begin{equation}\label{eq:split}
x^*-x_k=(x^*-x_{k+1})+(x_{k+1}-x_k).
\end{equation}
It suffices to prove that the two terms on the right-hand side are mutually orthogonal, and that the second does not vanish. For the first term, the consistency of \eqref{eq-1} and Lemma~\ref{lem:active} yield
\[
A_{\Lambda_{k+1}}(x^*-x_{k+1})=0.
\]
This implies that $(x^*-x_{k+1})\in\Null(A_{\Lambda_{k+1}})$. For the second term, \eqref{eq:xk} and $\Lambda_k\subset\Lambda_{k+1}$ give $x_k\in\Range(A_{\Lambda_k}^T)\subseteq\Range(A_{\Lambda_{k+1}}^T)$ and $x_{k+1}\in\Range(A_{\Lambda_{k+1}}^T)$, resulting in  $(x_{k+1}-x_k)\in\Range(A_{\Lambda_{k+1}}^T)$. From $\Range(A_{\Lambda_{k+1}}^T)\perp\Null(A_{\Lambda_{k+1}})$, it follows that
\begin{equation}\label{eq:pyth}
\|x^*-x_{k+1}\|_2^2=\|x^*-x_k\|_2^2-\|x_{k+1}-x_k\|_2^2.
\end{equation}
If Algorithm~\ref{alg:ocp} does not terminate at step $k$, then $r_k\neq 0$, and the greedy selection rule ensures that $r_k(i_{k+1})=b_{i_{k+1}}-a_{i_{k+1}}^Tx_k\neq 0$. In contrast, we observe from Lemma~\ref{lem:active} that $a_{i_{k+1}}^Tx_{k+1}=b_{i_{k+1}}$, that is, $x_{k+1}\neq x_k$. The assertion follows immediately.
\end{proof}

Actually, the result in Lemma~\ref{lem:monotone} can be further quantified. From Lemma~\ref{lem:active}, $a_{i_{k+1}}^T(x_{k+1}-x_k)=r_k(i_{k+1})$, and from $\|a_{i_{k+1}}\|_2=1$, we have $\|x_{k+1}-x_k\|_2\ge|r_k(i_{k+1})|=\|r_k\|_\infty$. Together this with \eqref{eq:pyth} gives
\[
\|e_{k+1}\|_2^2\le\|e_k\|_2^2-\|r_k\|_\infty^2.
\]
This implies that each iteration of Algorithm~\ref{alg:ocp} reduces the error at least as much as one such step. Moreover, it follows that $\|r_k\|_\infty^2\ge\|r_k\|_2^2/m$ and $\|r_k\|_2=\|Ae_k\|_2\ge\sigma_{\min}(A)\|e_k\|_2$, where $\sigma_{\min}(A)$ is the smallest nonzero singular value of $A$, and $e_k\in\Range(A^T)$. Therefore,
\[
\|e_{k+1}\|_2^2\le\Bigl(1-\frac{\sigma_{\min}^2(A)}{m}\Bigr)\|e_k\|_2^2 .
\]
Before termination, the error thus decreases at least linearly, with the factor that the same argument gives for the maximal-residual method.

\begin{theorem}\label{thm:finite}
Let  $\varepsilon=0$. Then Algorithm~\ref{alg:ocp} terminates in at most 
$\rank(A)$ iterations  with $r_k=0$, and its last iterate $x_k$ is the desired minimum-norm solution $x^*$ of \eqref{eq-1}.
\end{theorem}

\begin{proof}
Observe from Lemma~\ref{lem:indep}  that the rows of $A_{\Lambda_k}$ are linearly independent, and then $k\le\rank(A)$. Therefore, the first part of the assertion is true. If $\vert{}\Lambda_k\vert{}=\rank(A)$, then $\Range(A_{\Lambda_k}^T)=\Range(A^T)$. From Proposition~\ref{prop:omp}, we have $r_k=0$. The algorithm stops with $r_k=0$, or equivalently $Ax_k=b$. Since $x_k\in\Range(A^T)$, it follows from Proposition~\ref{prop:decomp} that $x_k=x^*$.
\end{proof}

\subsection{Computational aspects}\label{sec:impl}

Note that the core work of Algorithm~\ref{alg:ocp} is to update $y_{k+1}$. According to Proposition~\ref{prop:rowspace}, it is given by
\[
y_{k+1}=G_{k+1}^{-1}b_{\Lambda_{k+1}},
\]
where $G_{k+1}=A_{\Lambda_{k+1}}A_{\Lambda_{k+1}}^T$. To avoid the ${O}(k^3)$ operations of solving this system, there are two available alternatives.  The first one exploits $\Lambda_{k+1}=\Lambda_k\cup\{i_{k+1}\}$. We can recycle the previous $G_k^{-1}$ and apply the Schur complement of the positive definite matrix $G_{k+1}$. Let $\delta_{k+1}$ be the scalar defined by
\[
\delta_{k+1}=a_{i_{k+1}}^Ta_{i_{k+1}}-a_{i_{k+1}}^TA_{\Lambda_k}^Tt_{k+1},
\]
where $t_{k+1}=G_k^{-1}A_{\Lambda_k}a_{i_{k+1}}\in\mathbb{R}^{k}$. Given that $A_{\Lambda_k}^Tt_{k+1}$ is the orthogonal projection of $a_{i_{k+1}}$ onto $\Range(A_{\Lambda_k}^T)$, it follows that $\delta_{k+1}=\|a_{i_{k+1}}-A_{\Lambda_k}^Tt_{k+1}\|_2^2>0$. Then the $y_{k+1}$ can be computed recursively 
\begin{equation}\label{eq:schur}
\begin{aligned}
y_{k+1}&=\left(\begin{bmatrix}A_{\Lambda_k}\\ a_{i_{k+1}}^T\end{bmatrix}\begin{bmatrix}A_{\Lambda_k}^T& a_{i_{k+1}}\end{bmatrix}\right)^{-1}\begin{bmatrix}b_{\Lambda_k}\\ b_{i_{k+1}}\end{bmatrix}\\
&=\begin{bmatrix}G_k^{-1}+\frac{1}{\delta_{k+1}}t_{k+1}t_{k+1}^T & -\frac{1}{\delta_{k+1}}t_{k+1}\\[1mm] -\frac{1}{\delta_{k+1}}t_{k+1}^T & \frac{1}{\delta_{k+1}}\end{bmatrix}\begin{bmatrix}b_{\Lambda_k}\\ b_{i_{k+1}}\end{bmatrix}\\
&=\left(\begin{bmatrix}G_k^{-1}&0\\0&0\end{bmatrix}+\frac{1}{\delta_{k+1}}\begin{bmatrix}t_{k+1}\\-1\end{bmatrix}\begin{bmatrix}t_{k+1}^T&-1\end{bmatrix}\right)\begin{bmatrix}b_{\Lambda_k}\\ b_{i_{k+1}}\end{bmatrix}\\
&=\begin{bmatrix}y_k\\0\end{bmatrix}+\frac{1}{\delta_{k+1}}\begin{bmatrix}t_{k+1}\\-1\end{bmatrix}\begin{bmatrix}t_{k+1}^T&-1\end{bmatrix}\begin{bmatrix}b_{\Lambda_k}\\ b_{i_{k+1}}\end{bmatrix}\\
&=\begin{bmatrix}y_k\\0\end{bmatrix}+\frac{1}{\delta_{k+1}}\begin{bmatrix}t_{k+1}\\-1\end{bmatrix}\bigl(a_{i_{k+1}}^TA_{\Lambda_k}^TG_k^{-1}b_{\Lambda_k}-b_{i_{k+1}}\bigr)\\
&=\begin{bmatrix}y_k\\0\end{bmatrix}+\frac{1}{\delta_{k+1}}\begin{bmatrix}t_{k+1}\\-1\end{bmatrix}\bigl(a_{i_{k+1}}^Tx_k-b_{i_{k+1}}\bigr)\\
&=\begin{bmatrix}y_k\\0\end{bmatrix}-\frac{r_k(i_{k+1})}{\delta_{k+1}}\begin{bmatrix}t_{k+1}\\-1\end{bmatrix}.
\end{aligned}
\end{equation}
In addition to the residual update, one iteration requires only $O(nk+k^2)$ operations, and eliminates the need to re-factorize the matrix $G_{k+1}$.

The update scheme above, working with the matrices $G_k$ whose condition numbers are the squares of those of $A_{\Lambda_k}^T$, may be inefficient, especially in ill-conditioned cases.
Updating the iterates by performing a QR factorization of $A_{\Lambda_k}^T$ is more stable than the former.
Since $A_{\Lambda_k}^T$ is a full column rank matrix of size $n\times k$, it has a QR factorization $A_{\Lambda_k}^T=Q_kR_k$, where $Q_k=[q_1,q_2,\dots,q_k]\in\R^{n\times k}$ has orthonormal columns, i.e., $Q_k^TQ_k=I_k$, and $R_k$ is a nonsingular upper triangular matrix of order $k$. Substituting this into the explicit expression of $x_k$ gives
\begin{equation}\label{eq:xqr}
x_k=Q_kR_k(R_k^TQ_k^TQ_kR_k)^{-1}b_{\Lambda_k}=Q_kR_k^{-T}b_{\Lambda_k}.
\end{equation}
Let $z_k$ be the solution of the triangular system $R_k^Tz_k=b_{\Lambda_k}$. Then $x_k=Q_kz_k$. Since also $x_k=A_{\Lambda_k}^Ty_k=Q_kR_ky_k$, the $y_k$ can be obtained by the triangular system $R_ky_k=z_k$.  When the constraint $i_{k+1}$ is activated, one step of the Gram--Schmidt process gives
\begin{equation}\label{eq:gs}
\begin{aligned}
h_k&=Q_k^Ta_{i_{k+1}},\\
\tilde q_{k+1}&=a_{i_{k+1}}-Q_kh_k,\\
\rho_k&=\|\tilde q_{k+1}\|_2,\\
q_{k+1}&=\tilde q_{k+1}/\rho_k.
\end{aligned}
\end{equation}
It follows that $A_{\Lambda_{k+1}}^T=Q_{k+1}R_{k+1}$, where $Q_{k+1}=[Q_k,\ q_{k+1}]$ and $R_{k+1}=\begin{bmatrix}R_k&h_k\\0&\rho_k\end{bmatrix}$. Note that $\rho_k>0$ because $a_{i_{k+1}}$ is linearly independent of all rows of $A_{\Lambda_k}$. In exact arithmetic, the classical form \eqref{eq:gs} and the modified Gram--Schmidt process yield the same $h_k$ and $q_{k+1}$. From the first relation in \eqref{eq:gs}, we have
\[
h_k^Tz_k=(Q_k^Ta_{i_{k+1}})^Tz_k=a_{i_{k+1}}^T(Q_kz_k)=a_{i_{k+1}}^Tx_k.
\]
Then the last equation of the augmented triangular system $R_{k+1}^Tz_{k+1}=b_{\Lambda_{k+1}}$ determines the scalar $\zeta_{k+1}$, that is,
\begin{equation}\label{eq:zeta}
h_k^Tz_k+\rho_k\zeta_{k+1}=b_{i_{k+1}}\ \Longrightarrow\ \zeta_{k+1}=\frac{r_k(i_{k+1})}{\rho_k}.
\end{equation}
Putting these relations together, the updated iterate is of the form
\begin{equation}\label{eq:xupd}
\begin{aligned}
x_{k+1}&=Q_{k+1}z_{k+1}=\begin{bmatrix}Q_k&q_{k+1}\end{bmatrix}\begin{bmatrix}z_k\\ \zeta_{k+1}\end{bmatrix}\\
&=Q_kz_k+\zeta_{k+1}q_{k+1}=x_k+\frac{r_k(i_{k+1})}{\rho_k}\,q_{k+1}.
\end{aligned}
\end{equation}
Hence, $\|x_{k+1}-x_k\|_2=|r_k(i_{k+1})|/\rho_k$, which gives the decrease in \eqref{eq:pyth} explicitly.

\begin{remark}\label{rem:reorth}
In finite-precision arithmetic, the columns of 
$Q_k$ computed by \eqref{eq:gs} lose orthogonality through cancellation when $\rho_k\ll\|a_{i_{k+1}}\|_2$. To prevent this, we can utilize the Gram--Schmidt process with reorthogonalization \cite{Abdelmalek71}, for which a second reorthogonalization is enough in practice. If $\rho_k/\|a_{i_{k+1}}\|_2<\eta$ for a threshold $\eta$, e.g., $\eta\approx1/\sqrt2$, then we orthogonalize once more,
\begin{equation}\label{eq:reorth}
\tilde q_{k+1}\leftarrow\tilde q_{k+1}-Q_k(Q_k^T\tilde q_{k+1}),\qquad \rho_k\leftarrow\|\tilde q_{k+1}\|_2,
\end{equation}
and append the resulting coefficients to $h_k$. If $\rho_k$ remains smaller than a certain small multiple of the product of the unit round-off error and $\|a_{i_{k+1}}\|_2$, then the selected row is numerically dependent on the rows indexed by $\Lambda_k$. In this case, we terminate the iteration.  This prevents the algorithm from dividing by an extremely small value $\rho_k$. According to Lemma \ref{lem:indep}, under exact arithmetic, such a situation will not occur before termination. 
\end{remark}

From \eqref{eq:gs}--\eqref{eq:reorth}, we see that the $k$th iteration of Algorithm~\ref{alg:ocp} requires one matrix-vector product $Aq_{k+1}$ for the residual update $r_{k+1}=r_k-\zeta_{k+1}Aq_{k+1}$ ($2\,\nnz(A)$ flops), one Gram--Schmidt step ($4nk$ flops, or $8nk$ with reorthogonalization), and a search over the $m$ entries of $r_k$. After $s$ iterations, the total cost is about $2s\,\nnz(A)+2ns^2$ flops ($2s\,\nnz(A)+4ns^2$ with reorthogonalization at every step). The storage is $ns+s^2/2+m+n$ numbers for $Q_s$, $R_s$, $r_s$, and $x_s$. Both the cost and storage are governed by the number $s\le\rank(A)$ of selected rows; $m$ enters only through $\nnz(A)$ and the residual. One iteration of CGLS or LSQR costs about $4\,\nnz(A)$ flops and requires $O(m+n)$ storage. The number of iterations these methods need for a given accuracy depends on the distribution of the singular values of $A$. Therefore, Algorithm~\ref{alg:ocp} is competitive when $s$ is small and does not exceed the number of iterations 
 of those Krylov methods for the same accuracy.
 If $s$ is close to $\rank(A)$ and $\rank(A)$ is large, the cost approaches that of a QR factorization of $A_\Lambda^T$. Krylov subspace methods are preferable if only $x^*$ is required.

\subsection{Inequality constraints}\label{sec:ineq}

Recall that the solution  $x^*$ of the hard-margin problem \eqref{eq-3} is the minimum-norm solution of  $A_\Lambda x=\e_\Lambda$, where $\Lambda=\supp(y^*)$ consists of support vectors. Now we proceed as in Algorithm~\ref{alg:ocp}.
We select rows greedily and compute minimum-norm solutions of the selected subsystems. For an index set $W$, if $Ax_W\ge\e$ and $y_W\ge0$, then the minimum-norm solution $x_W=A_W^Ty_W$ of $A_Wx=\e_W$ equals $x^*$. The reason is that $x_W$ and the vector $y$ with $y_W$ on $W$ and zero elsewhere then satisfy \eqref{eq:kkt}. In this subsection, we let $A\in\R^{m\times(n+1)}$ be an arbitrary matrix, and its rows are not normalized, since scaling a row changes the constraint
$a_i^Tx\ge1$. For $a_i\neq0$, the constraint $a_i^Tx\ge1$ is a half-space, and, as usual for inequalities, it is called active at $x$ if $a_i^Tx=1$, that is, if $x$ lies on its bounding hyperplane.

To make the OCP method applicable to the hard-margin problem \eqref{eq-3}, two modifications are required. First, the constraint $i$ is violated only if $1-a_i^Tx>0$. A negative value means that the constraint holds with some slack. The selection takes the one-sided violation $(1-a_i^Tx)_+$ instead of $|r_k(i)|$. On the other hand, the coefficients $y_W$ must be nonnegative, so that a selected constraint may have to leave $W$ again. Here, the working set $W$ contains the indices of 
the constraints that are currently treated as equalities, that is, as hyperplanes. We organize both modifications around the invariant 
\begin{equation}\label{eq:inv}
x=A_W^Ty_W,\qquad y_W\ge0,\qquad A_Wx=\e_W,
\end{equation}
with $A_W$ of full row rank. It holds for $x=0$ and $W=\emptyset$. By \eqref{eq:inv}, $x$ is the minimum-norm solution $x_W$, and its coefficients are nonnegative, so only the first condition $Ax\ge\e$ may fail.

If $\max_{i\notin W}(1-a_i^Tx)\le0$, with the maximum over an empty set taken as $-\infty$, then $Ax\ge\e$, since $a_j^Tx=1$ for $j\in W$. In this case $x=x^*$, and the $y$ with $y_W$ on $W$ and zero elsewhere satisfies \eqref{eq:kkt}. Otherwise, we select the most violated constraint
\[
\ell=\argmax_{i\notin W}\,(1-a_i^Tx),\qquad \nu_\ell=1-a_\ell^Tx>0.
\]
In exact arithmetic, the restriction $i\notin W$ changes nothing, since $1-a_j^Tx=0$ for $j\in W$. In floating-point arithmetic, these values are rounding errors, which can be positive and even maximal. If there is no such restriction, then a constraint in $W$ could be selected again. The constraint $\ell$ enters with a coefficient $y_\ell$ that is initially zero, which leads to $x=A_W^Ty_W+y_\ell a_\ell$. Let
\[
a_\ell=A_W^Tu+z,\qquad A_Wz=0,\qquad u=(A_WA_W^T)^{-1}A_Wa_\ell,
\]
where $z$ is the part of $a_\ell$ that is orthogonal to the rows of $A_W$. We move along the path
\begin{equation}\label{eq:path}
x(\tau)=x+\tau z,\qquad y_W(\tau)=y_W-\tau u,\qquad y_\ell(\tau)=y_\ell+\tau,\qquad \tau\ge0 .
\end{equation}
From $a_\ell-A_W^Tu=z$, $A_Wz=0$ and $a_\ell^Tz=\|z\|_2^2$, it follows that
\[
x(\tau)=A_W^Ty_W(\tau)+y_\ell(\tau)a_\ell,\qquad A_Wx(\tau)=\e_W,\qquad 1-a_\ell^Tx(\tau)=\nu_\ell-\tau\|z\|_2^2 .
\]
Along the path, the constraints in $W$ stay active, and the violation of constraint $\ell$ decreases linearly if $z\neq0$. Only the coefficients $y_j$ with $u_j>0$ decrease, and that may limit the step. 

It should be pointed out that two events can end a step. If $z\neq0$, then the constraint $\ell$ becomes active at $\tau_f=\nu_\ell/\|z\|_2^2$. This shows that  $x(\tau_f)$ lies in  $\Range(A_{W\cup\{\ell\}}^T)$ and satisfies $A_{W\cup\{\ell\}}x=\e_{W\cup\{\ell\}}$. One can see that it is the minimum-norm solution of this system, i.e., the point that Step~C of Algorithm~\ref{alg:ocp} computes for the index set $W\cup\{\ell\}$. If $z=0$, then the violation does not decrease, and $\tau_f=\infty$. A coefficient in $W$ reaches zero at
\[
\tau_d=\min\{y_j/u_j\mid j\in W,\ u_j>0\},
\]
with $\tau_d=\infty$ if $u\le0$. We set $\tau=\min\{\tau_f,\tau_d\}$. If $\tau=\tau_f$, called a full step, then the index $\ell$ joins the working set; otherwise, called a partial step, the constraint $j^*$ that attains $\tau_d$ leaves $W$. Its coefficient is zero, so the relations $x=A_W^Ty_W+y_\ell a_\ell$ and $A_Wx=\e_W$ still hold for the smaller set. Since $\tau<\tau_f$, the constraint $\ell$ is still violated. We can repeat the step with the same $\ell$, with $\nu_\ell$ replaced by the current violation $1-a_\ell^Tx>0$. If $z=0$ and $u\le0$, neither event occurs, and \eqref{eq-3} is infeasible. For the sake of clarity, we list the four possible outcomes of a step in Table~\ref{tab:steps}.
\begin{table}[htb]
\caption{The four possible outcomes of a step of Algorithm~\ref{alg:ocpineq} for the selected constraint $\ell$.}\label{tab:steps}
\centering
\small
\begin{tabular}{@{}cccccp{5.0cm}@{}}
\toprule
Case & $z$ & $u$ & $\tau_f$ & $\tau_d$ & Action \\
\midrule
1 & $\neq0$ & any & $\nu_\ell/\|z\|_2^2$ & $\ge\tau_f$ &
Full step: $\ell$ joins $W$, and $x$ becomes the minimum-norm solution of $A_{W\cup\{\ell\}}x=\e_{W\cup\{\ell\}}$; next outer iteration. \\[2pt]
2 & $\neq0$ & some $u_j>0$ & $\nu_\ell/\|z\|_2^2$ & $<\tau_f$ &
Partial step: $y_{j^*}$ reaches zero and $j^*$ leaves $W$; $\ell$ is still violated, and the step is repeated with the same $\ell$. \\[2pt]
3 & $0$ & some $u_j>0$ & $\infty$ & finite &
Partial step with $x$ unchanged: weight moves from $y_W$ to $y_\ell$, and $j^*$ leaves $W$; the step is repeated with the same $\ell$. \\[2pt]
4 & $0$ & $u\le0$ & $\infty$ & $\infty$ &
Stop: no point satisfies the constraints in $W\cup\{\ell\}$, and \eqref{eq-3} is infeasible. \\
\bottomrule
\end{tabular}
\end{table}

The dual objective of \eqref{eq-3} given by
\[
D(y)=\e^Ty-\tfrac12\|A^Ty\|_2^2
\]
measures the progress. For every $y\ge0$, we have $D(y)\le\tfrac12\|x^*\|_2^2$, with equality at the solution, and the iteration preserves $y\ge0$. Using $x^Tz=y_\ell\|z\|_2^2$ and $u^T\e_W=u^TA_Wx=a_\ell^Tx-y_\ell\|z\|_2^2$, it follows that
\begin{equation}\label{eq:dualinc}
D(y(\tau))-D(y(0))=\tau \nu_\ell-\frac{1}{2}\tau^2\|z\|_2^2,
\end{equation}
which is positive for $0<\tau\le\tau_f$. When $\tau>0$, each step increases $D$. This is the basis of the termination proof. From the above results, we devise a variant of OCP for solving \eqref{eq-3}, see Algorithm \ref{alg:ocpineq}. The method can also be viewed as the dual active-set method of Goldfarb and Idnani \cite{GoldfarbIdnani83}, with the maximal-violation selection rule of Algorithm~\ref{alg:ocp}, applied to a strictly convex quadratic program with Hessian $I$. In exact arithmetic, we may take $\varepsilon=0$ in its stopping test, while in floating-point arithmetic, $\varepsilon$ is a small positive tolerance.

 \vspace{4mm}
\begin{breakablealgorithm}
\caption{OCP for $Ax\ge\e$}\label{alg:ocpineq}
\begin{algorithmic}[1]
\item[\algorithmicrequire] Input $A$ and the tolerance $\varepsilon\ge0$.
\State $x=0$, $W=\emptyset$
\Loop
\If{$\max_{i\notin W}(1-a_i^Tx)\le\varepsilon$} \textbf{stop} \EndIf
\State $\ell=\argmax_{i\notin W}\,(1-a_i^Tx)$,\quad $y_\ell=0$
\Loop
\State write $a_\ell=A_W^Tu+z$ with $A_Wz=0$
\State $\tau_f=(1-a_\ell^Tx)/\|z\|_2^2$ if $z\neq0$, and $\tau_f=\infty$ otherwise
\State $\tau_d=\min\{y_j/u_j\mid j\in W,\ u_j>0\}$, attained at $j^*$ ($\tau_d=\infty$ if $u\le0$)
\If{$\tau_f=\tau_d=\infty$} \textbf{stop}: \eqref{eq-3} is infeasible \EndIf
\State $\tau=\min\{\tau_f,\tau_d\}$,\quad $x\leftarrow x+\tau z$,\quad $y_W\leftarrow y_W-\tau u$,\quad $y_\ell\leftarrow y_\ell+\tau$
\If{$\tau=\tau_f$} $W\leftarrow W\cup\{\ell\}$; leave the inner loop \EndIf
\State $W\leftarrow W\setminus\{j^*\}$
\EndLoop
\EndLoop
\item[\algorithmicensure] Output $x$, $W$, and the multipliers $y_W$, or report that \eqref{eq-3} is infeasible.
\end{algorithmic}
\end{breakablealgorithm}
 \vspace{4mm}

\begin{proposition}\label{prop:ineq}
Let $\varepsilon=0$ and assume exact arithmetic.
Then Algorithm~\ref{alg:ocpineq} stops after finitely many iterations. 
If Algorithm~\ref{alg:ocpineq} stops because $\max_{i\notin W}(1-a_i^Tx)\le0$, then $x$ solves \eqref{eq-3}, and the vector $y$ with $y_W$ on $W$ and zeros elsewhere is an optimal multiplier vector; in particular, $W\subset\{i\mid a_i^Tx^*=1\}$ and $|W|\le\rank(A)\le n+1$. If it stops because $\tau_f=\tau_d=\infty$, then \eqref{eq-3} is infeasible.
\end{proposition}

\begin{proof}
First, we show that the invariant is maintained.
Within an inner loop, the relations $x=A_W^Ty_W+y_\ell a_\ell$, $y_W\ge0$, $y_\ell\ge0$, and $A_Wx=\e_W$ hold. The path \eqref{eq:path} keeps the first and the last terms
, the ratio test keeps $y_W\ge0$, and a removal keeps all of them since the removed coefficient is zero. A constraint is added only after a full step, where $z\neq0$, so that $a_\ell\notin\Range(A_W^T)$ and $A_{W\cup\{\ell\}}$ has full row rank. Hence \eqref{eq:inv} holds at the beginning of every outer iteration, and $|W|\le\rank(A)$.

If Algorithm \ref{alg:ocpineq} stops because $\max_{i\notin W}(1-a_i^Tx)\le0$, then $Ax\ge\e$, since $a_j^Tx=1$ for $j\in W$. By \eqref{eq:inv}, the vectors $x$ and $y$ satisfy \eqref{eq:kkt}.  From the fact that the constraints of \eqref{eq-3} are linear, the KKT conditions are sufficient for optimality, and $x=x^*$. If $z=0$ and $u\le0$, then $a_\ell=A_W^Tu$, and every $x'$ with $A_Wx'\ge\e_W$ satisfies $a_\ell^Tx'=u^TA_Wx'\le u^T\e_W=a_\ell^Tx<1$. This implies that no feasible point exists. 

It remains to prove finite termination. An inner loop consists of at most $|W|$ partial steps, since each of them removes an element of $W$. If $W$ becomes empty, then $z=a_\ell$, and the inner loop ends with a full step, or with the exit for infeasibility if $a_\ell=0$. It follows from \eqref{eq:dualinc} that the function $D$ does not decrease at any step, and
it increases strictly at a full step, where $\nu_\ell>0$ and $z\neq0$. After a full step, $x\in\Range(A_W^T)$ and $A_Wx=\e_W$. Therefore, $x$ is the minimum-norm solution $x_W$ of $A_Wx=\e_W$, and 
\[
D(y)=y_W^TA_Wx-\tfrac12\|x\|_2^2=\tfrac12\|x_W\|_2^2
\]
depends on $W$ only. The working sets accepted after full steps have strictly increasing values of  $\|x_W\|_2$. No working set is accepted twice, and the number of full steps is finite.
\end{proof}

\begin{remark}\label{rem:ineq}
If no constraint is ever removed, then every step is a full step. The iterates of Algorithm~\ref{alg:ocpineq} coincide with those of Algorithm~\ref{alg:ocp} applied to $Ax=\e$, with the one-sided violation in place of $|r_k(i)|$. The removal rule matters. Consider a simpler rule that appends the most violated constraint to $W$. Then we compute the coefficients from $(A_WA_W^T+\alpha I)y_W=\e_W$ with a small $\alpha>0$, and remove the constraint with the most negative coefficient until all coefficients are nonnegative. This rule neither keeps $A_W$ of full row rank nor needs to terminate. For example, given five points $p_1=(2.5,-0.6)$, $p_2=(-1.4,0.3)$, $p_3=(-2.5,1.3)$, $p_4=(-1.2,2.4)$, $p_5=(1.6,-2.6)$ with corresponding labels $q=(1,-1,-1,1,-1)$, for which $\rank(A)=3$, with $\alpha=10^{-10}$ and with ties in the selection broken by the smallest index. The working sets selected by the algorithm in sequence are $\{1\}$, $\{1,5\}$, $\{1,3,5\}$, and then in turn $\{3,4,5\}$, $\{2,4\}$, $\{1,2,4\}$, $\{3,4\}$, $\{3,4,5\}$, and so on. Throughout this iterative process, the corresponding dual objective function values $\frac12\|x_W\|_2^2$ are approximately $0.7008$, $0.6710$, $0.6820$, $0.6938$, and $0.7008$, respectively.  The same cycle is reached for every choice of the first constraint. The monotonicity used in the proof of Proposition~\ref{prop:ineq} fails, and after $\{3,4,5\}$ the rule adds the constraint $2$, and the four rows of $A_{\{2,3,4,5\}}$ are linearly dependent. On the same data and tie-breaking, Algorithm~\ref{alg:ocpineq} accepts $\{1\}$, $\{1,5\}$, and $\{1,3,5\}$, removes $1$ and accepts $\{3,4,5\}$, removes $5$ and accepts $\{2,3,4\}$, and stops, while $\frac12\|x_W\|_2^2$ increases from $0.0657$ to $0.7124$. These show that the iterations required by Algorithm~\ref{alg:ocpineq} are not bounded by $\rank(A)$, as they are for Algorithm~\ref{alg:ocp} by Theorem~\ref{thm:finite}. The recovery results of section~\ref{sec:recovery} do not apply to Algorithm~\ref{alg:ocpineq}. At termination, however, the KKT conditions certify the computed support, independently of the path.
\end{remark}

\section{Recovery of the support}\label{sec:recovery}
Notably, Theorem~\ref{thm:finite} bounds the number of iterations by $\rank(A)$. Assume that $x^*$ has a representation $x^*=A_S^Ty_S$ with $|S|=s\ll\rank(A)$. The theorem does not illustrate whether the set $\Lambda$ returned by Algorithm~\ref{alg:ocp} is then small, nor whether $\Lambda$ coincides with the support $S$ of this representation. From Proposition~\ref{prop:omp}, these are questions about OMP, and the answers known for OMP carry over to Algorithm~\ref{alg:ocp}.

The argument rests on the following observation. Since the rows of $A$ are normalized,  $|r_k(i)|=|\langle a_i,e_k\rangle|$ is the distance from $x_k$ to the $i$th hyperplane, and Algorithm~\ref{alg:ocp} selects the farthest hyperplane. It returns $S$ after $s$ iterations if and only if the index corresponding to the selected hyperplane belongs to $S$ at each iteration. As long as all selected rows belong to $S$, the run on $Ax=b$ makes the same selections as the run on the subsystem $A_Sx=b_S$. Moreover, by Proposition~\ref{prop:omp}, the iterates and errors produced by both runs
are the same. These errors lie in the $s$-dimensional subspace $\Range(A_S^T)$, and are independent of the rows outside of $S$.

We first turn this observation into a criterion for exact recovery. Under the exact-recovery conditions proposed by Tropp \cite{Tropp04}, this criterion holds for every coefficient vector $y_S$. Moreover, it also holds with high probability for rows uniformly distributed on the unit sphere. For the latter, the only ingredient that has to be adapted to the present setting is the lower-bound estimate for the smallest singular value $\sigma_s(A_S)$ of $A_S$.
This is because the rows of $A$ are normalized, so its entries are not independent.
The lower bound follows by applying diagonal scaling to the results obtained for the Gaussian case. For the sake of convenience, we let $\sigma_s(M)$ be the $s$th largest singular value of a matrix $M$.

\begin{lemma}\label{lem:criterion}
Let $S\subset\{1,\dots,m\}$, $|S|=s$, be such that $A_S$ has full row rank, let $x^*=A_S^Ty_S$ with $y_S\in\R^s$ having no zero entries, and let $b=Ax^*$. Then Algorithm~\ref{alg:ocp} with $\varepsilon=0$, applied to $A_Sx=b_S$, stops after exactly $s$ iterations at $x^*$.  Denote its errors by $e_0^*,\dots,e_{s-1}^*$. If
\begin{equation}\label{eq:criterion}
\max_{j\notin S}|\langle a_j,e_k^*\rangle|<\max_{i\in S}|\langle a_i,e_k^*\rangle|,\qquad k=0,\dots,s-1,
\end{equation}
then Algorithm~\ref{alg:ocp}, applied to $Ax=b$ with the same tie-breaking rule, makes the same selections as the run on $A_Sx=b_S$ and stops after exactly $s$ iterations with $\Lambda_s=S$ and $x_s=x^*$. Conversely, if \eqref{eq:criterion} holds for $k<k_0$, and the reverse strict inequality holds for $k=k_0\le s-1$, then the run on $Ax=b$ selects an index outside $S$ at iteration $k_0+1$.
\end{lemma}

\begin{proof}
From Theorem~\ref{thm:finite}, one can see that the run on  $A_Sx=b_S$ stops at $x^*$ after at most $s$ iterations. If it stopped after $k<s$ iterations, then $x^*$ is a combination of the $k$ rows indexed by $\Lambda_k\subsetneq S$. This contradicts the uniqueness of the coefficients of $x^*$ with respect to the linearly independent rows of $A_S$, since no entry of $y_S$ vanishes. So, this run stops after exactly $s$ iterations, and $e_k^*\neq0$ for $k<s$. Assume that the two runs have made the same selections in the first $k<s$ iterations. Then they have the same active set $\Lambda_k\subset S$ and iterate $P_{\Range(A_{\Lambda_k}^T)}x^*$. Therefore, the residual of the run on $Ax=b$ has the entries $\langle a_i,e_k^*\rangle$, $i=1,\dots,m$.

It follows from \eqref{eq:criterion} that the largest of their absolute values is attained only in $S$, and the run on $Ax=b$ selects the same index as the run on $A_Sx=b_S$. By induction, the two runs make the same $s$ selections. Since no index is selected again, $\Lambda_s=S$, and $x_s=x^*$. Under the counter-hypothesis, the two runs yield identical results for the first $k_0$ iterations, and at the  $k_0+1$th iterate, the largest absolute value of the residual entries 
occurs exclusively outside of $S$.
\end{proof}

Now we assume that the rows  $a_j$, $j\notin S$, are uniformly on $\Sph^{n-1}$ and independent of  $A_S$ and $y_S$. In this case, the probability that the equality in \eqref{eq:criterion} holds is zero, and hence, with probability one, \eqref{eq:criterion} is also a necessary condition for the run on $Ax=b$ to return $\Lambda_s=S$. The errors $e_k^*$ in \eqref{eq:criterion}, depending on $y_S$, are nonzero vectors in $\Range(A_S^T)$. Therefore, \eqref{eq:criterion} holds for every coefficient vector $y_S$ if $\max_{j\notin S}|\langle a_j,e\rangle|<\max_{i\in S}|\langle a_i,e\rangle|$ for every nonzero $e\in\Range(A_S^T)$.

\begin{theorem}\label{the:erc}
Let $S\subset\{1,\dots,m\}$, $|S|=s$, be such that $A_S$ has full row rank and
\begin{equation}\label{eq:erc}
\max_{j\notin S}\bigl\|(A_SA_S^T)^{-1}A_S\,a_j\bigr\|_1<1.
\end{equation}
Then for every $b=Ax^*$, with $x^*=A_S^Ty_S$ and $y_S\in\R^s$ having no zero entries, Algorithm~\ref{alg:ocp} with $\varepsilon=0$ stops after $s$ iterations with $\Lambda_s=S$ and $x_s=x^*$. The condition \eqref{eq:erc} holds whenever $s<\frac12(1+\mu^{-1})$, where $\mu=\max_{i\neq j}|a_i^Ta_j|$.
\end{theorem}

\begin{proof}
For $j\notin S$, one can deduce that the orthogonal projection of $a_j$ onto $\Range(A_S^T)$ is $A_S^Tc_j$, where $c_j=(A_SA_S^T)^{-1}A_Sa_j$. For every nonzero $e\in\Range(A_S^T)$,
it follows that
\[
\langle a_j,e
\rangle=\langle A_S^Tc_j,e
\rangle
=c_j^TA_Se,
\] and  $A_Se\neq0$. Putting this relation together with \eqref{eq:erc} gives \[
|\langle a_j,e\rangle|=|c_j^TA_Se|\le\|c_j\|_1\|A_Se\|_\infty<\max_{i\in S}|\langle a_i,e\rangle| .
\]
As shown in Lemma~\ref{lem:criterion}, the errors $e_k^*$ are nonzero and lie in $\Range(A_S^T)$. Then \eqref{eq:criterion} holds, which gives the first assertion.
For the last assertion, we let $G=A_SA_S^T=I+E$, where $I$ is an identity matrix of order $s$, the diagonal entries of $E$ are $0$, and its nondiagonal entries $E_{il}$ satisfy $|E_{il}|\leq \mu$.
Then $\|E\|_1\leq (s-1)\mu$
and $\|A_Sa_j\|_1=\sum_{i\in S}
a_i^Ta_j\leq s\mu$. If $(s-1)\mu< 1$, then $\|G^{-1}\|_1\leq \frac{1}{1-(s-1)\mu}$. Therefore, we have
$$\|c_j\|_1\leq\|G^{-1}\|_1\|A_Sa_j\|_1\leq \frac{s\mu}{1-(s-1)\mu}.$$
It is less than $1$ if and only if $(2s-1)\mu<1$. This condition also implies $(s - 1)\mu < 1$, which ensures that $G$ is invertible and $A_S$ has full row rank. 
\end{proof}

As seen, Theorem~\ref{the:erc} holds for every coefficient vector, but its coherence condition is restrictive. In particular, for $m>n$ unit vectors in $\R^n$, the Welch bound \cite{Welch74} gives $\mu\ge\sqrt{(m-n)/(n(m-1))}$, so that $s<\frac12(1+\mu^{-1})$ admits only $s$ of the order of $\sqrt n$. For rows uniformly distributed on the unit sphere, the criterion can hold with high probability for $s$ of the order of $n/\ln m$. The argument compares the two sides of \eqref{eq:criterion}. From a lower bound on $\sigma_s(A_S)$, some row in $S$ correlates with $e_k^*$ of at least the order of $\|e_k^*\|_2/\sqrt s$, whereas a row outside $S$ is independent of $e_k^*$, and its correlation with $e_k^*$ is typically of the order of $\|e_k^*\|_2/\sqrt n$. In what follows, we present a helpful lemma for the random recovery.

\begin{lemma}\label{lem:smin}
Let $A_S\in\R^{s\times n}$, $s<n$, have independent rows uniformly distributed on $\Sph^{n-1}$. Then, for every $t\ge0$,
\[
\PP\Bigl\{\sigma_s(A_S)\ge\frac{\sqrt n-\sqrt s-t}{\sqrt n+t}\Bigr\}\ \ge\ 1-(s+2)\,e^{-t^2/2}.
\]
In particular, if $n\ge16s$, then $\sigma_s(A_S)\ge2/5$ with probability at least $1-(s+2)e^{-n/32}$.
\end{lemma}

\begin{proof}
Let $g_1,\dots,g_s$ be independent  $N(0,I_n)$ vectors, $G_S=[g_1,\dots,g_s]^T$, and $D=\diag(\|g_1\|_2,\dots,\|g_s\|_2)$. Then $D^{-1}G_S$ has the same distribution as $A_S$. Since the assertion depends only on the distribution of $A_S$, we may assume that $A_S=D^{-1}G_S$. Then, for every $u\in\R^s$,
\[
\|A_S^Tu\|_2=\|G_S^TD^{-1}u\|_2\ \ge\ \sigma_s(G_S^T)\,\|D^{-1}u\|_2\ \ge\ \frac{\sigma_s(G_S^T)}{\max_i\|g_i\|_2}\,\|u\|_2.
\]
As seen from Corollary~5.35 in \cite{Vershynin12}, we have that $\sigma_s(G_S^T)\ge\sqrt n-\sqrt s-t$ with probability at least  $1-2e^{-t^2/2}$. Since $g\mapsto\|g\|_2$ is $1$-Lipschitz and $\EE\|g\|_2\le\sqrt n$, the Gaussian concentration 
\cite {Vershynin12} gives $\PP\{\|g_i\|_2>\sqrt n+t\}\le e^{-t^2/2}$ for each $i$, $i=1,\cdots, s$. A union bound over these $s+1$ events proves the first statement. If $n\ge16s$ and $t=\sqrt n/4$, the second follows immediately.
\end{proof}

The following result is the row-space form of the theorem of Tropp and Gilbert on OMP with admissible measurement matrices \cite[Theorem~6]{TroppGilbert07}. We give the argument in full because the constraints are explicit and the dictionary $A^T$ considered here has normalized columns rather than independent entries. 
In the random model, the support $S$ and the coefficients $y^*$ are commonly fixed first, and the rows of $A$ are drawn independently of them; otherwise, the results below are invalid.

\begin{theorem}\label{the:random}
Fix $S\subset\{1,\dots,m\}$ with $|S|=s$ and $y^*\in\R^m$ with $\supp(y^*)=S$. Let the rows of $A\in\R^{m\times n}$ be independent and uniformly distributed on   $\Sph^{n-1}$, independent of $y^*$, and set  $x^*=A^Ty^*$ and $b=Ax^*$. If $m\ge2$ and $n\ge50\,s\ln m$, then  with probability at least $1-\frac4m$,  Algorithm~\ref{alg:ocp} with $\varepsilon=0$ selects only indices
in $S$, stops after exactly $s$ iterations, and returns $\Lambda_s=S$, $x_s=x^*$, and $y_{\Lambda_s}=y_S^*$, up to the ordering of the entries.
\end{theorem}

\begin{proof}
Since $x^*\in\Range(A^T)$ and $Ax^*=b$, $x^*$ is the minimum-norm solution of \eqref{eq-1}. From Proposition~\ref{prop:omp}, for every realization of $A$, the event in question coincides with OMP applied to $\Phi=A^T$ and $x^*=\Phi y^*$, which recovers $y^*$ in $s$ iterations. Therefore,  following the argument in \cite{TroppGilbert07}, we show that \eqref{eq:criterion} holds with high probability. Let $e_0^*,\dots,e_{s-1}^*$ be the errors of Algorithm~\ref{alg:ocp} applied to $A_Sx=b_S$. They only depend on $A_S$ and $y^*$, and are independent of $\{a_j\}_{j\notin S}$. This independence is why we work with these errors,
rather than with those
from the run on $Ax=b$, which may depend on the rows outside $S$.

\emph{Rows in $S$.} From $e_k^*\in\Range(A_S^T)$, it follows that
\[
\max_{i\in S}|\langle a_i,e_k^*\rangle|\ \ge\ s^{-1/2}\|A_Se_k^*\|_2\ \ge\ s^{-1/2}\sigma_s(A_S)\,\|e_k^*\|_2.
\]
Let $\mathcal E=\{\sigma_s(A_S)\ge2/5\}$. If $n\ge50\,s\ln m $ and $m\ge 2$, then $n>16s$. Following by Lemma~\ref{lem:smin} and 
$s\ge1$
yields $\PP(\mathcal E^c)\le(s+2)e^{-n/32}\le(s+2)m^{-25s/16}\le3/m$. On $\mathcal E$, the submatrix $A_S$ has full row rank. It follows from Lemma~\ref{lem:criterion} that $\max_{i\in S}|\langle a_i,e_k^*\rangle|\ge\theta\|e_k^*\|_2$ for $k=0,\dots,s-1$, with $\theta=2/(5\sqrt s)$.

\emph{Rows outside $S$.} For a unit vector $u$ independent of $a_j$, the spherical cap estimate \cite[Lemma~2.2]{Ball97} yields $\PP\{|\langle a_j,u\rangle|\ge\theta\}\le2e^{-n\theta^2/2}$. Moreover, we have $n\ge\frac{25}{2}\,s\,(2\ln m+\ln2s)$ since $n\ge50\,s\ln m$ and $2s\le m^2$. With this relation, $\theta$, and a union bound over $j\notin S$, $k=0,\dots,s-1$, it follows that
\[
\PP\Bigl\{\max_k\max_{j\notin S}\frac{|\langle a_j,e_k^*\rangle|}{\|e_k^*\|_2}\ge\theta\ \Big|\ A_S\Bigr\}\ \le\ 2sm\,e^{-2n/(25s)}\ \le\ \frac1m.
\]
Consequently, on $\mathcal E$ and outside the event in the last display, \eqref{eq:criterion} holds, and Lemma~\ref{lem:criterion} gives $\Lambda_s=S$ and $x_s=x^*$. Since $A_S$ has full row rank, the coefficients of $x^*$ with respect to the rows of $A_S$ are unique, so $y_{\Lambda_s}=y_S^*$ up to the ordering of the entries. Overall, the probability of failure is at most $3/m+1/m=4/m$.
\end{proof}

\begin{remark}\label{rem:recovery}
 Notice that Theorems~\ref{the:erc} and \ref{the:random} present only sufficient conditions. The factor $\ln m$ in Theorem~\ref{the:random} comes from the union bound over the rows outside $S$. The constant $50$ in Theorem~\ref{the:random} has not been optimized. In numerical experiments, the recovery probability passes through $1/2$ at about $n\approx0.5\,s\ln m$, and these experiments do not demonstrate whether the factor $\ln m$ is necessary for Algorithm~\ref{alg:ocp}. For a given matrix and a set $S$, the condition \eqref{eq:erc} can be verified by performing a
factorization of $A_SA_S^T$ and computing the product of $A$ with $A_S^T$ of size $n\times s$. Both results pertain specifically to Algorithm~\ref{alg:ocp} involving equality constraints, while they are infeasible for Algorithm~\ref{alg:ocpineq}.
\end{remark}

\section{Numerical experiments}\label{sec:numerics}

In this section, we present numerical results on three kinds of problems to illustrate the identification of representing rows, and then show the effectiveness of our methods over many state-of-the-art Kaczmarz-type solvers for the consistent linear system \eqref{eq-1}.
Examples~\ref{ex:svm} and \ref{ex:SVMcale} are hard-margin SVM on real and synthetic data, in which the support vectors are the quantity of interest; Example~\ref{ex:random} illustrates Theorem~\ref{the:random} on random ensembles; and Example~\ref{ex:uf} uses sparse matrices from the SuiteSparse Matrix Collection, for which $x^*$ has no small representing set, to compare OCP with Kaczmarz-type methods and Krylov subspace methods. Denote by ``IT'' the number of iterations and by ``CPU'' the elapsed computing time in seconds. All experiments have been carried out in MATLAB R2026a on a personal computer with an Intel(R) Core(TM) i5-12400F CPU at 2.5\,GHz and 16\,GB memory.

\begin{example}\label{ex:svm}
We consider the hard-margin SVM \eqref{eq-3} on the banknote authentication data set \cite{UCI,Lohweg12} from the UCI Machine Learning Repository, which contains $1372$ samples with four features, including variance, skewness, kurtosis, and entropy of wavelet-transformed images, $610$ of them {in class $1$} ($q_i=1$) and $762$ {in class $0$} ($q_i=-1$). Each feature is standardized over the whole data set. We take the first $500$ samples of each class in the order of the data file, fit a linear soft-margin SVM with penalty parameter {$1000$}, and discard the samples whose functional margin with respect to this classifier does not exceed $0.1$. This leaves $m=987$ linearly separable samples, and problem \eqref{eq-3} has $n+1=5$ unknowns.

As stated in section~\ref{sec:intro}, this application requires a set $\Lambda$ of support vectors together with multipliers $y_\Lambda\ge0$ satisfying $x^*=A_\Lambda^Ty_\Lambda$, not only the separating vector $x^*\in\R^5$. We apply Algorithm~\ref{alg:ocpineq} and compare it with two solvers, namely the primal-dual interior-point method (PDIPM) of the function \texttt{quadprog} in MATLAB's Optimization Toolbox \cite{Mehrotra92,Wright97}, applied to the dual of the standard hard-margin SVM, and the sequential minimal optimization (SMO) method of the function \texttt{fitcsvm} in the Statistics and Machine Learning Toolbox \cite{FanChenLin05,ChangLin11}, with box constraint $10^5$. Both solvers treat the bias as a variable without regularization, so that their dual problems contain the equality constraint $\sum_iq_iy_i=0$, and SMO with a finite box constraint solves the soft-margin problem. These problems thus differ from \eqref{eq-3}; for the present data set, however, they have the same solution. We run each solver 20 times
after a warm-up call, and report their median time and other results in Table \ref{tab:svm}.

\begin{table}[htb]
\caption{Numerical results on Example~\ref{ex:svm}}\label{tab:svm}
\centering\footnotesize
\begin{tabular}{lccc}
\toprule
 & Algorithm~\ref{alg:ocpineq} & PDIPM & SMO \\
\midrule
problem solved & \eqref{eq-3} & standard SVM & standard SVM \\
$\|x-x^*\|_2$ & -- & $3.9\times10^{-6}$ & $6.8\times10^{-3}$ \\
multipliers $>10^{-5}$ & 5 & 10 & 5 \\
distinct support vectors & 5 & 5 & 5 \\
iterations & 11 additions, 6 removals & 13 & 11258 \\
CPU (s) & $2.9\times10^{-4}$& $1.4\times10^{-1}$ & $4.6\times10^{-2}$\\
\bottomrule
\end{tabular}
\end{table}

As seen from Table~\ref{tab:svm}, Algorithm~\ref{alg:ocpineq} stops after $11$ additions and $6$ removals, and returns
\[
x^*=[-9.2101,\ -14.2877,\ -11.3946,\ -1.7167,\ -5.2092]^T.
\]
The working set never exceeds $5$ constraints, and at termination $W$ consists of five constraints with multipliers $146.59$, $75.25$, $17.98$, $178.14$, and $30.93$. The KKT conditions \eqref{eq:kkt}, which are sufficient for optimality, confirm that this vector solves \eqref{eq-3}. They hold to working precision, with $\min_i(a_i^Tx^*-1)=-1.5\cdot10^{-14}$ and $\|x^*-A_W^Ty_W\|_2=1.1\cdot10^{-13}$. In particular, the classifier separates all $987$ samples. We observe that nine constraints are active at $x^*$, yet they correspond to only five distinct samples. The reason is that
 one support vector occurs twice and another one four times in the present data. Therefore, the optimal multiplier vector in \eqref{eq-3} is not unique since the KKT conditions \eqref{eq:kkt} hold for every nonnegative distribution of the multipliers among identical rows. The rows of the five distinct support vectors are linearly independent. This implies that every optimal multiplier vector has at least five nonzero entries, and Algorithm~\ref{alg:ocpineq} returns one with exactly five.

 Indeed, under PDIPM and SMO, $x^*$ is also the solution to the standard hard-margin SVM; that is, it is a vertex of the feasible set. In this case, the KKT conditions of the standard problem hold at $x^*$ with multipliers $139.29$, $71.59$, $16.57$, $165.20$, and $29.10$ on the five distinct support vectors, which satisfy $\sum_iq_iy_i=0$. Thus, we take $x^*$ as the reference solution for both solvers, and the symbol ``--'' in Table~\ref{tab:svm} means that the distance $\|x-x^*\|_2$ does not apply to Algorithm~\ref{alg:ocpineq} itself. The SMO solver returns five multipliers above the threshold $10^{-5}$, on the same five support vectors, while PDIPM returns ten, with nine on the nine active constraints, with equal values on identical rows, and one of size $10^{-4}$ on an inactive constraint. Additionally, we can see that the computing time required by Algorithm~\ref{alg:ocpineq} is much less than that of the other solvers, which is favorable for such problems.
\end{example}

\begin{example}\label{ex:SVMcale}
To examine how Algorithm~\ref{alg:ocpineq} depends on the number of samples, we consider the hard-margin SVM \eqref{eq-3} on linearly separable synthetic data. The samples are $p_i=g_i+\gamma q_iw_0$, $i=1,\dots,m$, where the vectors $g_i\sim N(0,I_n)$ are independent, $w_0\in\R^n$ is a fixed unit vector, $q_i=\operatorname{sign}(w_0^Tg_i)$, and $\gamma=0.1$. The slab can separate the two classes, defined by $|w_0^Tp|<\gamma$. The coefficient matrix, with rows $a_i^T=q_i[p_i^T,1]$, has full column rank $n+1\ll m$. Since the points $p_i$ have a continuous distribution, with probability one at most $n+1$ constraints are active at $x^*$. We take $n=10$ and $n=50$, and $m=10^4,10^5,10^6$.

Here, we compare Algorithm~\ref{alg:ocpineq} with \texttt{quadprog}, that is, the interior-point-convex method, now applied to \eqref{eq-3} itself with its $n+1$ variables and $m$ inequality constraints, and with \texttt{fitcsvm}, that is, SMO, standard SVM with free bias, with the same settings as in Example~\ref{ex:svm}.  Since \texttt{quadprog} takes about nine minutes for $m=10^4$, we run it for this size only, and we run SMO for the next value of $m$ only if the previous run took less than one minute. The results are recorded in Table~\ref{tab:SVMcale}.

\begin{table}[htb]
\caption{Example~\ref{ex:SVMcale}: hard-margin SVM on $m$ separable samples in $\R^n$. For Algorithm~\ref{alg:ocpineq}: number of support vectors $|W|$ at termination, number of constraints added and removed. For \texttt{quadprog} applied to \eqref{eq-3} and for \texttt{fitcsvm}: number of multipliers above $10^{-5}$, and distance to the solution $x^*$ of \eqref{eq-3} computed by Algorithm~\ref{alg:ocpineq}; $\dagger$ marks a run stopped at the iteration limit $10^6$ of SMO without convergence, and -- denotes that a solver was not run for that size.}\label{tab:SVMcale}
\centering\footnotesize\setlength{\tabcolsep}{3.5pt}
\begin{tabular}{rrrrrrrrrrrr}
\toprule
 & & \multicolumn{4}{c}{Algorithm~\ref{alg:ocpineq}} & \multicolumn{3}{c}{\texttt{quadprog}} & \multicolumn{3}{c}{\texttt{fitcsvm}} \\
\cmidrule(lr){3-6}\cmidrule(lr){7-9}\cmidrule(lr){10-12}
$n$ & $m$ & $|W|$ & added & removed & CPU & CPU & mult. & $\|x-x^*\|_2$ & CPU & SV & $\|x-x^*\|_2$ \\
\midrule
10 & $10^{4}$ & 11 & 44 & 33 & 0.003 & 566.5 & 11 & 9.5e-12 & 4.5 & 12 & 7.8e-03 \\
10 & $10^{5}$ & 11 & 58 & 47 & 0.027 & -- & -- & -- & 135.8 & 16 & 3.5e-04 \\
10 & $10^{6}$ & 11 & 64 & 53 & 0.715 & -- & -- & -- & -- & -- & -- \\
50 & $10^{4}$ & 51 & 228 & 177 & 0.056 & 543.3 & 51 & 3.1e-11 & 36.6 & 52 & 4.9e-03 \\
50 & $10^{5}$ & 51 & 340 & 289 & 0.821 & -- & -- & -- & 584.4$^\dagger$ & 60 & 2.9e-03 \\
50 & $10^{6}$ & 51 & 413 & 362 & 15.7 & -- & -- & -- & -- & -- & -- \\
\bottomrule
\end{tabular}
\end{table}

We see from Table~\ref{tab:SVMcale} that Algorithm~\ref{alg:ocpineq} stops with $|W|=n+1$ in all six cases, as expected, with the KKT conditions satisfied to working precision. The number of iterations grows slowly with $m$, from $77$ to $117$ for $n=10$ and from $405$ to $775$ for $n=50$, that is, $7$ to $15$ times $n+1$. Since the residual costs $O(mn)$ flops per iteration, one would expect the time to grow in proportion to $m$ times the number of iterations. However, from $m=10^5$ to $m=10^6$,  it grows by factors of about $26$ with $n=10$ and $19$ with $n=50$, while the number of iterations grows only by factors of about $1.1$ and $1.2$. We believe that this is a memory effect. For $m=10^6$, the matrix $A$ occupies $88$ MB with $n=10$ and $408$ MB with $n=50$, and each residual evaluation reads it once. For $m=10^6$, the time is $0.7$ seconds for $n=10$ and $16$ seconds for $n=50$.

The interior-point solver finds the same solution and set of positive multipliers in $14$ and $13$ iterations with $\|x-x^*\|_2\le4\cdot10^{-11}$, but requires $566$ and $543$ seconds for $m=10^4$, since each iteration works with all $m$ constraints. One can see that SMO converges in $5\cdot10^4$ to $3.4\cdot10^5$ iterations and $4.5$ to $135.8$ seconds on three of the problems and is stopped at its iteration limit $10^6$ after $584.4$ seconds on the fourth. Its solutions differ from $x^*$ by $4\cdot10^{-4}$ to $8\cdot10^{-3}$ and violate $Ax\ge\e$ by $3\cdot10^{-4}$ to $8\cdot10^{-4}$. The margin $1/\|w\|_2$ agrees with that of $x^*$ to four digits, but its set of support vectors has one to nine more elements than $W$ and differs from $W$ in all four cases. For both values of $n$, Algorithm~\ref{alg:ocpineq} requires less time for $m=10^6$ than the other two solvers for $m=10^4$. The working set never exceeds $n+1$ constraints, and $m$ enters the cost only through one pass over $A$ per iteration.
\end{example}

\begin{figure}[htb]
\centering
\includegraphics[width=0.48\textwidth]{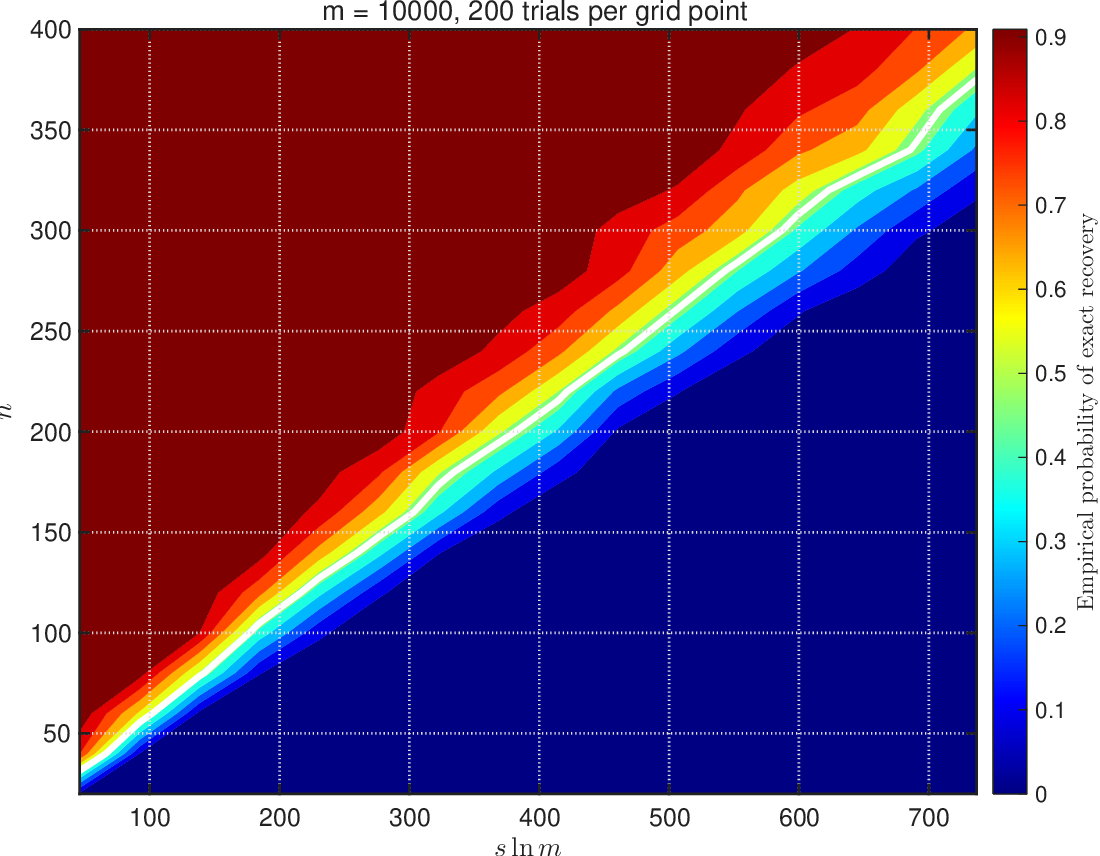}
\ \ \includegraphics[width=0.48\textwidth]{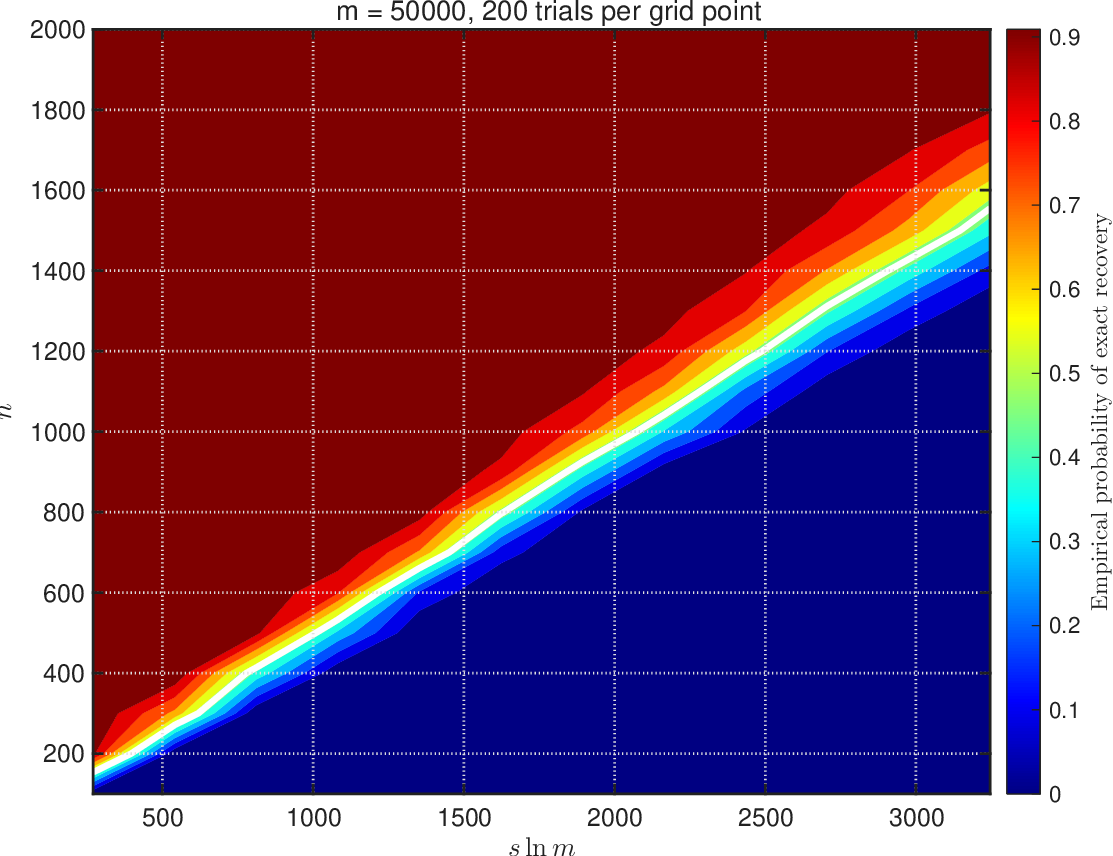}
\caption{Empirical probability of exact recovery in Example~\ref{ex:random} with $m=10000$ (Left), and $m=50000$ (Right), where
the white curve is the level $1/2$.}\label{fig:phase1}
\end{figure}

\begin{example}\label{ex:random}
In this example, we draw an $m\times n$ standard Gaussian random matrix and normalize its rows to unit $\ell_2$-norm, so that they are independent and uniformly distributed on $\Sph^{n-1}$. A support set $S\subset\{1,\dots,m\}$ with $|S|=s$ is sampled uniformly at random, $y_S^*\sim N(0,I_s)$ is drawn, and $x^*=A_S^Ty_S^*$, $b=Ax^*$. This is the random model of Theorem~\ref{the:random}. At each point $(s,n)$ of a grid, the empirical recovery rate of Algorithm~\ref{alg:ocp} is estimated from independent trials; a trial is successful if the algorithm selects exactly the rows in $S$ in $s$ iterations, in which case $x_s=x^*$ up to rounding. For $s\le n$, Lemma~\ref{lem:criterion} and the paragraph following it show that, up to an event of probability zero, a trial is successful if and only if $\max_{j\notin S}|\langle a_j,e_k^*\rangle|<\max_{i\in S}|\langle a_i,e_k^*\rangle|$ for $k=0,\dots,s-1$, where $e_k^*$ are the errors of the run on $A_Sx=b_S$; for $s>n$, no trial is successful, since Algorithm~\ref{alg:ocp} stops after at most $n$ iterations. Since $e_k^*\in\Range(A_S^T)$, only the projections of the rows $a_j$, $j\notin S$, onto this $s$-dimensional subspace are needed, and conditionally on $A_S$ these projections are distributed as the first $s$ coordinates of a uniform point on $\Sph^{n-1}$. We use this equivalent form to evaluate each trial at a cost of $O(ms^2+ns^2)$.

\begin{table}[htb]
\caption{Example~\ref{ex:random} with $m=10000$: empirical probability of exact recovery at the grid points $n$ around the transition, written as $n$ (probability), and the smallest grid value $n_{1/2}$ with probability above $1/2$. Each entry is based on 3000 trials, or on 500 trials when all of them had the same outcome.}\label{tab:phase1}
\centering\footnotesize
\setlength{\tabcolsep}{4pt}
\begin{tabular}{rrrccccc}
\toprule
$s$ & $s\ln m$ & $n_{1/2}$ & \multicolumn{5}{c}{$n$ (probability)} \\
\midrule
10 & 92 & 60 & 20 (0.000) & 40 (0.080) & 60 (0.652) & 80 (0.937) & 100 (0.992) \\
20 & 184 & 120 & 80 (0.081) & 100 (0.456) & 120 (0.781) & 140 (0.936) & 160 (0.981) \\
40 & 368 & 200 & 160 (0.102) & 180 (0.310) & 200 (0.595) & 220 (0.783) & 240 (0.895) \\
60 & 553 & 300 & 260 (0.269) & 280 (0.471) & 300 (0.639) & 320 (0.794) & 340 (0.882) \\
80 & 737 & 380 & 320 (0.122) & 340 (0.252) & 360 (0.418) & 380 (0.553) & 400 (0.690) \\
\bottomrule
\end{tabular}
\end{table}

We take $m=10000$ with $n=20,40,\cdots,400$ and $m=50000$ with $n=100,200,\cdots,2000$, and report the numerical results in Tables \ref{tab:phase1} and \ref{tab:phase2}. For $m=10000$ and $s=5,10,\cdots,80$, and for $m=50000$ with $s=25,50,\cdots,300$, Figure~\ref{fig:phase1} illustrates the empirical recovery probability over the whole grid, estimated from $200$ trials per grid point, with the horizontal axis parameterized by $\kappa=s\ln m$ and the white curve marking the level $1/2$. For five values of $s$ in each case, Tables~\ref{tab:phase1} and \ref{tab:phase2} list more accurate estimates at the grid points around the transition, along with the smallest grid value $n_{1/2}$ at which the estimated probability exceeds $1/2$. Since $n_{1/2}$ is a grid value, it may overestimate the point where the probability passes $1/2$ by up to one grid step, that is, by $20$ for $m=10000$ and $100$ for $m=50000$.

We see that, for every $s$, the recovery probability rises from $0$ to $1$ over a range of $n$ comparable to $n_{1/2}$ itself. For $m=10000$ and $s=10$ from $0.080$ at $n=40$ to $0.992$ at $n=100$, for $s=40$ from $0.102$ at $n=160$ to $0.895$ at $n=240$, and for $s=80$ it is still $0.690$ at $n=400$, the largest dimension of the grid. The values $n_{1/2}$ grow roughly linearly in $s\ln m$. The ratio $n_{1/2}/(s\ln m)$ lies between $0.49$ and $0.74$, largest for the smallest $s$, where the grid step is a sizable fraction of $n_{1/2}$, and a least-squares fit through the origin gives $n_{1/2}\approx0.53\,s\ln m$ for $m=10000$ and $0.51\,s\ln m$ for $m=50000$. For $m=50000$ and $s=300$, for instance, the probability is zero for all $n\le1200$, $0.450$ at $n=1500$, and $0.865$ at $n=1800$, and $n_{1/2}=1600$ is close to the fitted value $0.51\,s\ln m\approx1650$. This is consistent with the dependence on $s$ and $m$ in Theorem~\ref{the:random}.
\end{example}
\begin{table}[htb]
\caption{Example~\ref{ex:random} with $m=50000$: empirical probability of exact recovery at the grid points $n$ around the transition, written as $n$ (probability), and the smallest grid value $n_{1/2}$ with probability above $1/2$. Each entry is based on 1000 trials, or on 200 trials when all of them had the same outcome.}\label{tab:phase2}
\centering\footnotesize
\setlength{\tabcolsep}{4pt}
\begin{tabular}{rrrccccc}
\toprule
$s$ & $s\ln m$ & $n_{1/2}$ & \multicolumn{5}{c}{$n$ (probability)} \\
\midrule
25 & 270 & 200 & 100 (0.009) & 200 (0.954) & 300 (1.000) & 400 (1.000) & 500 (1.000) \\
75 & 811 & 500 & 300 (0.004) & 400 (0.464) & 500 (0.927) & 600 (1.000) & 700 (1.000) \\
125 & 1352 & 700 & 500 (0.003) & 600 (0.194) & 700 (0.689) & 800 (0.925) & 900 (0.988) \\
175 & 1893 & 1000 & 800 (0.082) & 900 (0.456) & 1000 (0.767) & 1100 (0.916) & 1200 (0.971) \\
300 & 3246 & 1600 & 1400 (0.150) & 1500 (0.450) & 1600 (0.650) & 1700 (0.802) & 1800 (0.865) \\
\bottomrule
\end{tabular}
\end{table}

\begin{example}\label{ex:uf}
Consider solving \eqref{eq-1} with coefficient matrices from the SuiteSparse Matrix Collection \cite{DavisHu11}. 
These benchmark matrices arise mainly in combinatorial optimization and computer vision. We remove the zero rows in \texttt{relat7b} and \texttt{rel8} and normalize the remaining rows of each matrix. Table~\ref{tab:uf} gives the dimensions after this removal. The matrices \texttt{relat7b}, \texttt{rel8} and \texttt{specular} are rank deficient, and \texttt{specular} is also ill-conditioned. We generate a sparse vector $y$ by the MATLAB function \texttt{sprandn} with density $0.1$ and set $x^*=A^Ty$ and $b=Ax^*$, so that $x^*\in\Range(A^T)$ is the minimum-norm solution for the rank-deficient matrices as well. Then $x^*$ is a combination of about $m/10$ rows, more than $\rank(A)$, and no representation by few rows is to be expected. All iterative solvers use the same stopping criterion $\|x_k-x^*\|_2<10^{-6}$. The Kaczmarz-type methods and OCP are stopped after $1000$ seconds, and the Krylov subspace methods after $50n$ iterations. The symbol $\dagger$ indicates that a run was stopped by one of these limits without meeting the criterion. The compared methods are as follows:
\begin{itemize}
\item[$\star$] MR: the maximal-residual method \eqref{eq:iter}--\eqref{eq:select} of Agmon and Motzkin--Schoenberg, i.e., the memoryless counterpart of OCP, with the residual recomputed per iteration;
\item[$\star$] RaBK: randomized average block Kaczmarz method with adaptive step size and uniform block sampling \cite{Necoara19}, the rows being partitioned into $k_b=\lceil\|A\|_2^2\rceil$ consecutive blocks;
\item[$\star$] FGBK: fast greedy block Kaczmarz method \cite{XiaoYinZheng23} with $p=2$ and relaxation parameter $0.05$;
\item[$\star$]CGNE \cite{Craig55}, LSQR \cite{PaigeSaunders82}, and CGLS \cite{Bjorck96}, in their standard forms, without preconditioning and reorthogonalization;
\item[$\star$] OCP: Algorithm~\ref{alg:ocp} implemented with \eqref{eq:gs}--\eqref{eq:reorth}.
\end{itemize}
For OCP, we also report the fraction $|\Lambda|/n$ of rows selected at termination.

\begin{table}[htb]
\caption{Example~\ref{ex:uf}: iterations (IT), CPU time in seconds, and the error $\|x-x^*\|_2$ (EN); $\dagger$ marks runs stopped by the time limit of $1000$s or, for the Krylov methods, by the limit of $50n$ iterations. For OCP, $|\Lambda|/n$ is the fraction of rows selected.}\label{tab:uf}
\centering\footnotesize\setlength{\tabcolsep}{1.5pt}
\begin{tabular}{llccccc}
\toprule
Matrix & & \texttt{abtaha2} & \texttt{relat7b} & \texttt{rel8} & \texttt{mesh\_deform} & \texttt{specular} \\
$m\times n$ & & $37932\times331$ & $20500\times1045$ & $206040\times12347$ & $234023\times9393$ & $477976\times1600$ \\
$\nnz(A)$ & & 137228 & 81355 & 821839 & 853829 & 7647040 \\
\midrule
MR & IT & 1801 & 60046 & 298479$^\dagger$ & 333785$^\dagger$ & 78996$^\dagger$ \\
 & CPU & 0.52 & 10.20 & 1000.00 & 1000.00 & 1000.00 \\
 & EN & 9.9e-07 & 9.9e-07 & 1.5e-03 & 4.6e-05& 4.0e-01 \\
 RaBK & IT & 6294 & 9651 & 37711 & 32111 & 16344147$^\dagger$ \\
 & CPU & 0.17 & 0.40 & 17.47 & 18.26 & 1000.00 \\
 & EN & 9.8e-07 & 9.8e-07 & 9.7e-07 & 9.9e-07 & 6.6e-02 \\
FGBK & IT & 176 & 566 & 1762 & 7623 & 32117$^\dagger$ \\
 & CPU & 0.18 & 0.34 & 15.99 & 73.78 & 1000.00 \\
 & EN & 9.6e-07 & 9.7e-07 & 9.9e-07 & 9.8e-07 & 7.9e-02 \\
CGLS & IT & 40 & 76 & 117 & 223 & 80000$^\dagger$ \\
 & CPU & 0.02 & 0.02 & 0.38 & 0.66 & 1489.27 \\
 & EN & 5.7e-07 & 8.5e-07 & 9.2e-07 & 9.9e-07 & 1.2e-05 \\
LSQR & IT & 40 & 76 & 117 & 223 & 80000$^\dagger$ \\
 & CPU & 0.02 & 0.02 & 0.36 & 0.64 & 1538.55 \\
 & EN & 5.7e-07 & 8.5e-07 & 9.3e-07 & 9.9e-07 & 1.2e-05 \\
CGNE & IT & 39 & 74 & 115 & 218 & 80000$^\dagger$ \\
 & CPU & 0.02 & 0.02 & 0.37 & 0.56 & 1486.22 \\
 & EN & 9.6e-07 & 9.6e-07 & 9.6e-07 & 9.9e-07 & 9.4e-06 \\
OCP & IT & 323 & 995 & 6679$^\dagger$ & 8071$^\dagger$ & 1422 \\
 & CPU & 0.17 & 2.41 & 1000.00 & 1000.00 & 24.92 \\
 & EN & 7.4e-07 & 6.7e-07 & 4.8e+00 & 8.7e-01& 9.1e-07 \\
 & $|\Lambda|/n$ & 0.98 & 0.95 & 0.54 & 0.86 & 0.89 \\
\bottomrule
\end{tabular}
\end{table}

We implement all testing methods for solving \eqref{eq-1}, and list
the numerical results in Table~\ref{tab:uf}. As seen, for \texttt{specular}, the ill-conditioned matrix, Algorithm~\ref{alg:ocp} reaches the error criterion after $1422$ iterations and $24.92$ seconds, and it is the only method in this table that does so. The MR, RaBK, and FGBK are stopped by the time limit with errors between $6.6\cdot10^{-2}$ and $4.0\cdot10^{-1}$. The errors of the Krylov subspace methods stagnate at about $10^{-5}$, and the tolerance is not reached within $50n$ iterations. In exact arithmetic, OCP and the three Krylov methods would all terminate after at most $\rank(A)$ iterations, so that the difference is caused by rounding errors. Unlike the Krylov methods, OCP follows Remark~\ref{rem:reorth} for reorthogonalization, and selects $1422$ row in this matrix. Since the matrix has $158$ zeros columns, that is, $\rank(A)\le1442$, the QR factorization of $A_\Lambda^T$ in OCP almost covers the whole row space of $A$. To show the effect of reorthogonalization further,  we also ran LSQR with one-sided reorthogonalization, in which every new right vector of the Golub--Kahan bidiagonalization is orthogonalized twice against all previous ones. With the same $x^*$, $b$, and stopping criterion, this variant yields the  criterion after $1415$ iterations and $34$ seconds, with the error $8.5\times10^{-7}$, whereas LSQR without reorthogonalization has the error $2.2\times10^{-2}$ after the same number of iterations. This means that the result on \texttt{specular} reflects the reorthogonalization in OCP, rather than an advantage of OCP as a solver. Precondtioning, which we have not tried, is another manner to speedup the Krylov methods. For a rank-deficient $A$, however, a right preconditioner $M$ leads to the solution that minimizes $\|Mx\|_2$ instead of $x^*$. Note also that the final error $9.1\times10^{-7}$ of OCP is only slightly below the tolerance, so the outcome of this comparison depends on the choice of the tolerance.

Actually, Algorithm~\ref{alg:ocp} also reaches the error criterion on \texttt{abtaha2} and \texttt{relat7b}. On these three matrices, the number of selected rows is $89\%$ (\texttt{specular}), $95\%$ (\texttt{relat7b}), and $98\%$ (\texttt{abtaha2}) of $n$, and the number of iterations is much smaller than those of MR and RaBK. The comparison with MR, the memoryless counterpart of OCP, shows the gain from projecting onto the intersection of all selected constraints. On \texttt{abtaha2} and \texttt{relat7b}, OCP requires about $6$ and $60$ times fewer iterations than MR and about one-third and one-quarter of its computing time. However, FGBK requires fewer iterations than OCP on \texttt{abtaha2} and \texttt{relat7b}, and the block methods RaBK and FGBK are as fast as OCP on \texttt{abtaha2} and faster on \texttt{relat7b}.
On the \texttt{rel8} and \texttt{mesh\_deform} matrices, which have the largest $n$, OCP is stopped by the time limit after $6679$ and $8071$ iterations, and MR is stopped by this limit as well, with smaller errors, while RaBK and FGBK meet the error criterion at about $16$ to $74$ seconds. The Krylov subspace methods, which also terminate in at most $\rank(A)$ iterations in exact arithmetic, need $39$ to $223$ iterations on the first four matrices and are faster there than OCP and the Kaczmarz-type methods, by about one order of magnitude on \texttt{abtaha2} and by more than one order of magnitude on the other three matrices.

We see from this example that, on the three matrices on which it converges, Algorithm~\ref{alg:ocp} is faster than its memoryless counterpart MR.  Of course, $x^*$ has no small representing set here, so that $|\Lambda|$ is a sizable fraction of $n$, and the cost approaches that of a QR factorization of $A_\Lambda^T$, as noted in section~\ref{sec:impl}. If only $x^*$ is required, the Krylov subspace methods are cheaper on the first four matrices. None of the other methods, however, provides the set $\Lambda$ together with the coefficients $y_\Lambda$. The identification of this set when $x^*$ has a representing set of $s\ll m$ rows is illustrated by Examples~\ref{ex:svm}--\ref{ex:random}.
\end{example}

\section{Conclusion}\label{sec:conclusion}

In this paper, we proposed an orthogonal constraint pursuit method for identifying a small set of rows of a consistent linear system whose span
contains the minimum-norm solution $x^*$, along with the coefficients of that representation. This greedy method selects the constraint with the maximum violation, as in Motzkin's relaxation method, and projects the iterate onto the intersection of all constraints selected so far. We proved that the proposed method is equivalent to OMP applied to the dictionary $A^T$ and the signal $x^*$, and that it terminates in at most $\rank(A)$ iterations.  We also presented the probabilistic guarantee that a representation with $s$ rows, whose support and coefficients are chosen independently of $A$, is identified in exactly $s$ iterations when $n\ge Cs\ln m$, and the deterministic exact recovery condition.

The method is intended for systems whose rank is far smaller than the number of rows, in which $x^*$ is a combination of at most $\rank(A)$ rows and the identity of these rows, rather than the solution alone, is required. For $s$ iterations, the cost is $O(s\,\nnz(A)+ns^2)$ flops and the storage is $O(ns+m)$, and the selected rows come with the coefficients of the representation. When $x^*$ is represented by only a large number of rows, as in the sparse matrices of Example~\ref{ex:uf}, whose ranks are in the hundreds or thousands, the cost approaches that of a QR factorization, and Krylov subspace methods are preferable for computing $x^*$.

For inequality constraints, we gave a variant with constraint removal, the dual active-set method of Goldfarb and Idnani with the maximal-violation selection rule. It terminates finitely, and its output is certified by the KKT conditions. The recovery results do not cover this variant, and conditions under which its number of iterations can be bounded in terms of the number of support vectors remain to be studied.

\medskip
\noindent\textbf{Contributions:} All authors contributed equally to Conceptualization, Formal analysis, Investigation, Methodology, Software, Validation, Writing, Review, Editing. All authors have read and agreed to the published version of the manuscript.

\medskip
\noindent\textbf{Data Availability Statement:} The banknote authentication data set and the test matrices are publicly available from the UCI Machine Learning Repository \cite{UCI,Lohweg12} and the SuiteSparse Matrix Collection \cite{DavisHu11}.

\medskip
\noindent\textbf{Conflict of Interest:} The authors declare no conflict of interest.

\medskip
\noindent\textbf{Declaration of generative AI:}
During the preparation of this manuscript, the authors used Claude (Anthropic) to improve the language and readability of the text. 
The authors assume responsibility for all content.


\begin{thebibliography}{99}

\bibitem{SeberLee} G.A.F. Seber and A.J. Lee, \emph{Linear Regression Analysis}, 2nd ed., John Wiley \& Sons, Hoboken, NJ, 2003.

\bibitem{BDM13} C. Boutsidis, P. Drineas, and M. Magdon-Ismail, \emph{Near-optimal coresets for least-squares regression}, IEEE Trans. Inform. Theory, 59 (2013), pp. 6880--6892.

\bibitem{Natarajan95} B.K. Natarajan, \emph{Sparse approximate solutions to linear systems}, SIAM J. Comput., 24 (1995), pp. 227--234.



\bibitem{CandesTao05} E.J. Cand\`es and T. Tao, \emph{Decoding by linear programming}, IEEE Trans. Inform. Theory, 51 (2005), pp. 4203--4215.

\bibitem{RudelsonVershynin05} M. Rudelson and R. Vershynin, \emph{Geometric approach to error-correcting codes and reconstruction of signals}, Int. Math. Res. Not., 2005 (2005), pp. 4019--4041.


\bibitem{Tropp04} J.A. Tropp, \emph{Greed is good: Algorithmic results for sparse approximation}, IEEE Trans. Inform. Theory, 50 (2004), pp. 2231--2242.


\bibitem{LustigDonohoSantos} M. Lustig, D.L. Donoho, and J.M. Santos, et al., \emph{Compressed sensing MRI}, IEEE Signal Process. Mag., 25 (2008), pp. 72--82.

\bibitem{BruntonBruntonProctor} B.W. Brunton, S.L. Brunton, and J.L. Proctor, et al.,  \emph{Sparse sensor placement optimization for classification}, SIAM J. Appl. Math., 76 (2016), pp. 2099--2122.

\bibitem{TsangKwokCristianini} I.W. Tsang, J.T. Kwok, and P.M. Cheung, \emph{Core vector machines: Fast SVM training on very large data sets}, J. Mach. Learn. Res., 6 (2005), pp. 363--392.

\bibitem{Brooks11} J.P. Brooks, \emph{Support vector machines with the ramp loss and the hard margin loss}, Oper. Res., 59 (2011), pp. 467--479.

\bibitem{MangasarianMusicant99} O.L. Mangasarian and D.R. Musicant, \emph{Successive overrelaxation for support vector machines}, IEEE Trans. Neural Networks, 10 (1999), pp. 1032--1037.

\bibitem{StrohmerVershynin09} T. Strohmer and R. Vershynin, \emph{A randomized Kaczmarz algorithm with exponential convergence}, J. Fourier Anal. Appl., 15 (2009), pp. 262--278.

\bibitem{NeedellTropp14} D. Needell and J.A. Tropp, \emph{Paved with good intentions: analysis of a randomized block Kaczmarz method}, Linear Algebra Appl., 441 (2014), pp. 199--221.

\bibitem{Censor81} Y. Censor, \emph{Row-action methods for huge and sparse systems and their applications}, SIAM Rev., 23 (1981), pp. 444--466.

\bibitem{Necoara19} I. Necoara, \emph{Faster randomized block Kaczmarz algorithms}, SIAM J. Matrix Anal. Appl., 40 (2019), pp. 1425--1452.

\bibitem{MiaoWu22} C.Q. Miao and W.T. Wu, \emph{On greedy randomized average block Kaczmarz method for solving large linear systems}, J. Comput. Appl. Math., 413 (2022), 114372.

\bibitem{BaiWu18} Z.Z. Bai and W.T. Wu, \emph{On greedy randomized Kaczmarz method for solving large sparse linear systems}, SIAM J. Sci. Comput., 40 (2018), pp. A592--A606.

\bibitem{DuHayamiZheng21} Y.S. Du, K. Hayami, and N. Zheng, et al., \emph{Kaczmarz-type inner-iteration preconditioned flexible GMRES methods for consistent linear systems}, SIAM J. Sci. Comput., 43 (2021), pp. S345--S366.

\bibitem{JiangWuJiang23} Y.T. Jiang, G. Wu, and L. Jiang, \emph{A semi-randomized Kaczmarz method with simple random sampling for large-scale linear systems}, Adv. Comput. Math., 49 (2023), 20.

\bibitem{DNR25} M. Derezi\'nski, D. Needell, and E. Rebrova, et al., \emph{Randomized Kaczmarz methods with beyond-Krylov convergence}, SIAM J. Matrix Anal. Appl., 46 (2025), pp. 2558--2588.

\bibitem{JeongNeedellWu26} H. Jeong and D. Needell, \emph{A federated Kaczmarz algorithm}, BIT Numer. Math., 66 (2026), 43.

\bibitem{Craig55} E.J. Craig, \emph{The $N$-step iteration procedures}, J. Math. Phys., 34 (1955), pp. 64--73.



\bibitem{PaigeSaunders82} C.C. Paige and M.A. Saunders, \emph{LSQR: An algorithm for sparse linear equations and sparse least squares}, ACM Trans. Math. Software, 8 (1982), pp. 43--71.

\bibitem{Bjorck96} \AA. Bj\"orck, \emph{Numerical Methods for Least Squares Problems}, 2nd ed., SIAM, Philadelphia, PA, 2024.

\bibitem{PRK93} Y.C. Pati, R. Rezaiifar, and P.S. Krishnaprasad, \emph{Orthogonal matching pursuit: Recursive function approximation with applications to wavelet decomposition}, in Proc. 27th Asilomar Conf. Signals, Systems and Computers, 1993, pp. 40--44.

\bibitem{DMA97} G. Davis, S. Mallat, and M. Avellaneda, \emph{Adaptive greedy approximations}, Constr. Approx., 13 (1997), pp. 57--98.

\bibitem{TroppGilbert07} J.A. Tropp and A.C. Gilbert, \emph{Signal recovery from random measurements via orthogonal matching pursuit}, IEEE Trans. Inform. Theory, 53 (2007), pp. 4655--4666.

\bibitem{NeedellVershynin09} D. Needell and R. Vershynin, \emph{Uniform uncertainty principle and signal recovery via regularized orthogonal matching pursuit}, Found. Comput. Math., 9 (2009), pp. 317--334.

\bibitem{NeedellTropp09} D. Needell and J.A. Tropp, \emph{CoSaMP: Iterative signal recovery from incomplete and inaccurate samples}, Appl. Comput. Harmon. Anal., 26 (2009), pp. 301--321.



\bibitem{Agmon54} S. Agmon, \emph{The relaxation method for linear inequalities}, Canad. J. Math., 6 (1954), pp. 382--392.

\bibitem{MotzkinSchoenberg54} T.S. Motzkin and I.J. Schoenberg, \emph{The relaxation method for linear inequalities}, Canad. J. Math., 6 (1954), pp. 393--404.



\bibitem{GoldfarbIdnani83} D. Goldfarb and A. Idnani, \emph{A numerically stable dual method for solving strictly convex quadratic programs}, Math. Programming, 27 (1983), pp. 1--33.



\bibitem{Saad03} Y. Saad, \emph{Iterative Methods for Sparse Linear Systems}, 2nd ed., SIAM, Philadelphia, PA, 2003.

\bibitem{DeLoeraHaddockNeedell17} J.A. De Loera, J. Haddock, and D. Needell, \emph{A sampling Kaczmarz--Motzkin algorithm for linear feasibility}, SIAM J. Sci. Comput., 39 (2017), pp. S66--S87.

\bibitem{HaddockNeedell19} J. Haddock and D. Needell, \emph{On Motzkin's method for inconsistent linear systems}, BIT Numer. Math., 59 (2019), pp. 387--401.

\bibitem{RZE08} R. Rubinstein, M. Zibulevsky, and M. Elad, \emph{Efficient implementation of the K-SVD algorithm using batch orthogonal matching pursuit}, Cs Technion,  40 (2008), pp. 1-15.


\bibitem{Abdelmalek71} N.N. Abdelmalek, \emph{Round off error analysis for Gram--Schmidt method and solution of linear least squares problems}, BIT Numer. Math., 11 (1971), pp. 345--367.

\bibitem{Welch74} L.R. Welch, \emph{Lower bounds on the maximum cross correlation of signals}, IEEE Trans. Inform. Theory, 20 (1974), pp. 397--399.

\bibitem{Vershynin12} R. Vershynin, \emph{Introduction to the non-asymptotic analysis of random matrices}, in Compressed Sensing: Theory and Applications, Y.C. Eldar and G. Kutyniok, eds., Cambridge University Press, Cambridge, 2012, pp. 210--268.

\bibitem{Ball97} K. Ball, \emph{An elementary introduction to modern convex geometry}, Flavors of Geometry, 31 (1997), pp. 1--58.

\bibitem{XiaoYinZheng23} A.Q. Xiao, J.F. Yin, and N. Zheng, \emph{On fast greedy block Kaczmarz methods for solving large consistent linear systems}, Comput. Appl. Math., 42 (2023), 119.

\bibitem{UCI} A. Asuncion and D. Newman, \emph{UCI Machine Learning Repository}, University of California, Irvine, 2007.

\bibitem{Lohweg12} V. Lohweg, \emph{Banknote Authentication}, UCI Machine Learning Repository, 2012,  https://doi.org/10.24432/C55P57.

\bibitem{Mehrotra92} S. Mehrotra, \emph{On the implementation of a primal-dual interior point method}, SIAM J. Optim., 2 (1992), pp. 575--601.

\bibitem{Wright97} S.J. Wright, \emph{Primal-dual interior-point methods}, SIAM, Philadelphia, PA, 1997.

\bibitem{FanChenLin05} R.E. Fan, P.H. Chen, and C.J. Lin, \emph{Working set selection using second order information for training support vector machines}, J. Mach. Learn. Res., 6 (2005), pp. 1889--1918.

\bibitem{ChangLin11} C.C. Chang and C.J. Lin, \emph{LIBSVM: A library for support vector machines}, ACM Trans. Intell. Syst. Technol., 2 (2011), pp. 1--27.

\bibitem{DavisHu11} T.A. Davis and Y. Hu, \emph{The University of Florida sparse matrix collection}, ACM Trans. Math. Software, 38 (2011), pp. 1--25.


\end{thebibliography}
\end{document}